\documentclass[fleqn,11pt]{article}

\usepackage{color}
\usepackage{amsmath}
\usepackage{amsfonts}
    {\par \noindent {\bf Proof:}}%
    {\par \indent}
    {\par \noindent {\bf Counterexample:}}%
    {\par \indent}

    {\par \noindent {\bf Remark: }}%
    {\par \indent}

\numberwithin{equation}{section}

\newcommand{\be}{\begin{equation}}
\newcommand{\bes}{\begin{displaymath}}
\newcommand{\ee}{\end{equation}}
\newcommand{\ees}{\end{displaymath}}

\usepackage{latexsym}

\usepackage{rotate,epsfig}

\def\E {{\rm E}}

\begin{document}

\vspace{-4cm}
\title{{Failure Inference and Optimization for Step Stress Model Based on Bivariate Wiener Model}}
\author{{ S. Shemehsavar}$^{\dag *}$\footnote{$^*$ Corresponding author;  \newline
E-mail addresses: {\it shemehsavar@khayam.ut.ac.ir} (S.
Shemehsavar), {\it mortez.amini@ut.ac.ir} (Morteza Amini).} ~~~and
{Morteza Amini }$^{\dag}$\\ \mbox{}\\
{$^{\dag}${\small Department of Statistics, School of
Mathematics, Statistics and Computer Sciences,}}\\
{{\small College of Science, University of Tehran, P.O. Box 14155-6455, Tehran,
Iran}}\vspace{2mm}\\} \maketitle
\begin{abstract}
In this paper, we consider the situation under a life test, in which
the failure time of the test units are not related deterministically
to an observable stochastic time varying covariate. In such a case, the joint distribution of failure
time and a marker value would be useful for modeling the step
stress life test. The problem of
accelerating such an experiment is considered as the main aim of this paper. We present a step stress accelerated model based on a bivariate
Wiener process with one component as the latent
(unobservable) degradation process, which determines the failure
times and the other as a marker process, the
degradation values of which are recorded at times of failure.
Parametric inference based on the proposed model is discussed and
the optimization procedure for obtaining the optimal time for
changing the stress level is presented. The optimization criterion
is to minimize the approximate variance of the maximum likelihood
estimator of a percentile of the products' lifetime distribution.
\end{abstract}
\noindent {\bf Keywords}:  Bivariate normal, Fisher information matrix,
Inverse Gaussian distribution.\\
\textbf{AMS subject classification: 62N05, 60K10}

\section{ Introduction  }

The lifetime experiments have received attention recently, partly
because the high reliability of the manufactured products is
important in the current intense economical competition between
trading firms. Over time, several lifetime tests for assessing the
lifetime probability distribution of the products are developed,
ranging from simple Constant Stress Life Test (CSLT) to the Step
Stress Accelerated Degradation Test (SSADT). Two useful survey of available results are given in the books
of Nelson ,1990 and Bagdonavicius and Nikulin, 2010. For some recent papers concerning the lifetime experiments see
Pan and Balakrishnan, 2010, Pan et al., 2011, Jin, 2011, Simino et al., 2012 and Wang et al., 2012.

Life tests usually deal
with models for which failure occurs when an observable
degradation process crosses a threshold level. However, there are
practical situations in which the failure time of the test units are
not related deterministically to an observable marker covariate. In such a case, the
joint distribution of the failure time and a marker process would be
useful for modeling the step stress life test. Joint models for
marker evolution and failure are proposed in the literature under
the simple constant stress life tests, including
Jewell and Kalbfleisch, 1996, who examine jump processes for markers and an additive
relationship between the marker and the failure time hazard function, and
Yashin and Manton, 1997, who consider diffusion processes for
markers along with a quadratic relationship between the hazard
function and markers.

In most cases, the information about the latent
(unobservable) degradation path can only be obtained
using the related marker(s) and the fact that a failure occurs when
the latent degradation process crosses a known threshold. Whitmore
et al., 1998, proposed a constant stress bivariate Wiener model
in which one component represents the marker and the second, which
is latent, determines the failure time.

The constant stress life tests are usually very costly, since they
require destroying a considerable number of products for testing at
each level of stress. To handle this problem, Step Stress
Accelerated Life Tests (SSALT) were proposed as an economic
alternative to the constant stress life tests. In a SSALT framework,
each product is first tested, subject to a pre-determined stress
level for a specified duration, and the failure data are collected.
A product which survived until the end of the first step was again
tested at a higher stress level and for a different time duration.
The experiment is repeated for a specified number of stress levels
and terminated at a pre-determined censoring time. The
constant stress bivariate Wiener model proposed by Whitmore et al., 1998
is as well a costly experiment. Although
censoring in this model decreases the total time of the experiment
it does not solve the problem of efficiency. To handle this
problem, we consider a SSALT design under the
bivariate Wiener model.

An essential problem in an SSALT design is to determine
the optimal time for changing the stress level by the experimenter.
The problem of optimizing the test design have been extensively
studied in recent years. Three commonly used optimization criteria
are the minimum Approximated variance (Avar) of the Maximum
Likelihood Estimators (MLE) of reliability, Mean Time To Failure
(MTTF) and the quantiles of the population. For surveys of recent
results in optimization of life test designs, see in particular Tang et al., 2004, Liao and Tseng, 2006
and Tseng et al., 2009.

In this paper, we reconstruct the model proposed by Whitmore et al., 1998 in a SSALT framework. Such a generalized model is clearly
more economic than the constant stress model of Whitmore et al., 1998, since constant stress experiment requires destroying a
considerable number of products at each level of stress.
The Maximum likelihood and Bayesian estimation of the
parameters of the proposed model are discussed.
Next, we determine the optimal stress changing time by
minimizing the Avar of the MLE of the 100$p^{\mbox{th}}$ percentile
of the products' life time distribution.

The rest of this paper is organized as follows. In
Section 2, we introduce the SSALT model with a bivariate Wiener
process and derive the joint distribution of failure times and the
marker process. Parametric Inference based on the proposed model is
discussed in Section 3. The optimization criterion is described in
Section 4. Finally an illustrative example is
presented in Section 5.

\section{The Model}

Consider a two-dimensional Wiener diffusion process $\{(X(r),
Y(r))\}$, for $r\geq0$ with $(X(0), Y(0))=(0,0)$ (see Cox and Miller, 1965). In other words, under the normal stress
level $S_0$
\[(X(r),Y(r))|S_0\sim N_2(r{\mu_X}_0,r{\mu_Y}_0,r\sigma^2_X,r\sigma^2_Y,\rho),\]
where $N_2$ stands for the bivariate normal distribution. Assume
further that ${\mu_X}_0\geq0$, which guarantees the degradation
process $X(r)$ to be stochastically increasing in $r$.

The component $X(r)$ assumed to be a degradation process that
represents the level of deterioration of an item. An item fails
as soon as $X(r)$ reaches a threshold $D>0$.
This first passage time of the degradation process through the threshold is
denoted by a random variable $T$, namely
\begin{equation}\label{t}
T=\inf\{t|X{(t)}\geq D\}.
\end{equation}
The failure time $T$ follows an inverse Gaussian distribution (see for instance, Chhikara and
Folks, 1989), with the cumulative distribution
function (cdf) under the normal stress level $S_0$ as follows
\begin{equation}\label{gt}
G_0(t)=\Phi\left(\sqrt{\frac{1}{\sigma_X^2t}}({\mu_X}_0t-D)\right) +
exp\left\{\frac{2{\mu_X}_0D}{\sigma_X^2}\right\}
\Phi\left(-\sqrt{\frac{1}{\sigma_X^2t}}({\mu_X}_0t+D)\right),
\end{equation}
where $\Phi$ is the cdf of the standard normal distribution.

The degradation process $X(r)$ is assumed to be
unobservable. The component $Y(r)$ represents a marker process that
is correlated with the degradation process and tracks its progress.
Thus, results of the experiment are based on
observations on the marker process, supplemented by failure times of
failed items. We focus on the situation where marker measurements
are taken only at the failure or censoring times.

Consider the above bivariate process to model a SSALT problem. Under
a SSALT, each item is first tested subject to a stress level $S_1$
($S_1>S_0$) for a specified duration
$[0,\tau_1)$. If the item does not fail, it is tested again at a
higher stress level $S_2$ ($S_2>S_1$) for another
specified duration $[\tau_1,\tau_2)$. The experiment is continued
until the time $C$, under $m\geq 2$ stress levels
$S_m>S_{m-1}>\cdots>S_2>S_1$.
 The stress level of the experiment is then defined as
$$S =\left\{ \begin{array}{l l} S_1 &\mbox{for}\;\; 0\leq t<\tau_1\\
 S_2  &\mbox{for}\;\; \tau_1\leq t<\tau_2\\
 \vdots& \\
 S_m  &\mbox{for}\;\; \tau_{m-1}\leq t<C,
\end{array}\right.
$$
where the pre-specified values
$0<\tau_1<\tau_2<\cdots<\tau_{m-1}<C$ are called the {\it stress
changing times}.
\input{epsf}
\epsfxsize=4in \epsfysize=3in
\begin{figure}[!hbtp]
\centerline{\epsffile{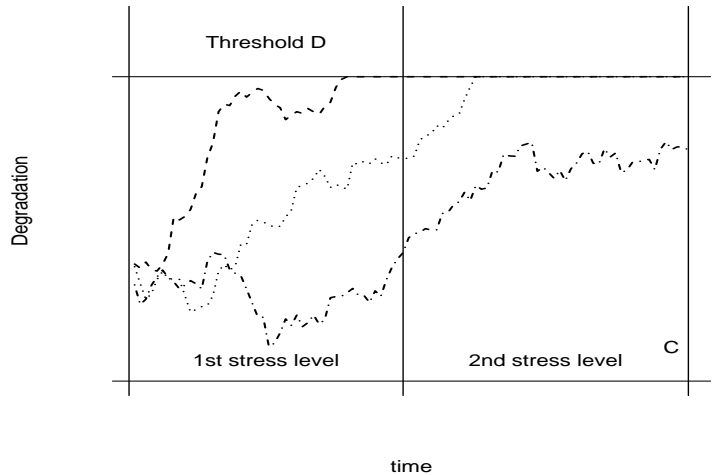}}
 \caption{\scriptsize A sampler degradation path. Three different paths are showed: a failed item under the first stress level (dashed line),
 a failed item under the second stress level (dotted line) and a survived (censored) item (dash-dotted line).}
\end{figure}

Under a SSALT model, each item has two possible observation
outcomes during the period $(0,C]$:
\begin{itemize}
\item \textbf{Surviving (Censored) item:} The item survives
to the censoring time $C$ at which a marker level of $Y(C)=y(C)$ is
recorded. This occurrence constitutes a censored observation of
failure time with $T>C$.
\item \textbf{Failing item:}
The items fails at some time $T=t$ during the period $(0,C]$ and a
marker level of $Y(T)=y(t)$ is recorded at the moment of failure.
\end{itemize}
\subsection{The distribution of failure time and marker covariate}
For the aforementioned plan, under the stress $S_j,$
for $j=1,2,\ldots,m$, we have
$$(X(r),Y(r))|S_j\sim
N_2(r\mu_{X_{j}},r\mu_{Y_{j}},r\sigma^2_X,r\sigma^2_Y,\rho).$$

Assume further that the Arrhenius reaction model is used to model
the relationship between the location parameters $\mu_{X_{j}}$ and
 $\mu_{Y_{j}}$  and the temperature stress $S_j$, that is
 \begin{equation}\label{mu}
\mu_{X_{j}}=\exp\left(a+\frac{b}{273+S_j}\right),\;\mu_{Y_{j}}=\exp\left(c+\frac{d}{273+S_j}\right),\;j=0,1,\ldots,m.
\end{equation}
 Consider any sample path of the
component $X$, under stress $S_j$, over a
 time interval $(0,r]$ and partition this sample path at
 arbitrary time points $0=r_0<r_1\cdots<r_k=r,\;k\geq 1$. Let
$\Delta r_i=r_i-r_{i-1}$ and $\Delta x_i=x(r_i)-x(r_{i-1})$, for
$i=1,\cdots,k$. Denote the set of realized increments $\{\Delta
x_1,\ldots,\Delta x_k\}$ by $P$. Then we have clearly
\begin{equation}\label{normal}
 Y{(r)}|P\sim N(\mu_{y.x(r)},r\sigma^2_Y(1-\rho^2)),
\end{equation}
where for $j=1,\ldots,m,$
\[\mu_{y.x(r)} =\mu_{j}(y,r)+\rho\frac{\sigma_{Y}}{\sigma_{X}}[x(r)-\mu_{j}(x,t)],\quad \tau_{j-1}\leq r<\tau_j,\]
\[\mu_{j}(y,t)=\mu_{Y_{j}}(t-\tau_{j-1})+\sum_{k=1}^{j-1}\mu_{Y_{k}}(\tau_k-\tau_{k-1}),\]
\[\mu_{j}(x,t)=\mu_{X_{j}}(t-\tau_{j-1})+\sum_{k=1}^{j-1}\mu_{X_{k}}(\tau_k-\tau_{k-1}),\]
$\tau_0=0$ and $\tau_m=C$.

The conditional distribution in \eqref{normal} is the same as the
conditional distribution $Y{(r)}|x(r)$. Hence \eqref{normal}
holds for any sample path of $X$.

Therefore, for a surviving path, the conditional
distribution of the marker given the degradation at the censoring
time $C$ is as follows
$$Y{(C)}|X(C)\sim N(\mu_{y.x(C)},C\sigma^2_Y(1-\rho^2)).$$
The resulting conditional probability density function
(p.d.f) of the surviving path then is
\begin{equation}\label{p1}
p_{1}(y|x;\theta)=C^{-1/2}\sigma^{-1}_Y(1-\rho^2)^{-1/2}\phi\left(C^{-1/2}\sigma^{-1}_Y(1-\rho^2)^{-1/2}(y-\mu_{y.x(C)})\right),
\end{equation}
 where $\phi$ is the pdf of the standard normal distribution and
$$\theta=(\mu_{X1},\ldots,\mu_{Xm},\mu_{Y1},\ldots,\mu_{Ym}
,\sigma_X^2,\sigma_Y^2,\rho).$$

For a failing item at time $t$, the distribution of
$Y{(t)}|x(t)$ is equal to \eqref{normal} with $r$ replaced by $t$
and $x(r)$ replaced by $x(t)=D$. The corresponding p.d.f
then is

\[p_{2}(y|t;\theta)=t^{-1/2}\sigma^{-1}_Y(1-\rho^2)^{-1/2}\phi\left(t^{-1/2}\sigma^{-1}_Y(1-\rho^2)^{-1/2}\left(y-\mu_{j}(y,t)+\rho\frac{\sigma_{Y}}{\sigma_{X}}[D-\mu_{j}(x,t)]\right)\right).\]

A similar argument to that in Lu, 1995 can be used to derive the p.d.f. of a
surviving item, that is ${\rm P}(X(C)=x,T>C)$, as follows
\begin{equation}\label{p3}
p_{3}(x)=\frac{1-exp(-\frac{2D(D-x)}{\sigma^2_{X}C})}{\sigma_{X}\sqrt{C}}
\phi(\frac{x-\mu_m(x,C)}{\sigma_{X}\sqrt{C}}),\quad-\infty<x<D.
\end{equation}
It is easy to verify that the p.d.f. of $T$ in
\eqref{t} of a failing item is
\[f_{T}(t|T<C) =\sum_{j=1}^{m}\frac{D}{\sqrt{(2\pi\sigma^2_{X}t^3)}}exp(-\frac{(D-\mu_{X_{j}}t)^2}{2\sigma^2_{X}t})I_{(\tau_{j-1},\tau_j)}(t),\]
where
$$I_A(t)=\left\{\begin{array}{lr}1\;\;\mbox{if}&t\in A\\0\;\;\mbox{if}&t\notin A,\end{array}\right.$$
$\tau_0=0$ and $\tau_m=C$.

We combine the preceding results to obtain the p.d.f. for each
type of observation outcome, as follows:

For a censored item which survives beyond time $C$, the joint p.d.f.
of the marker $Y(C)$ and the latent degradation $X(C)$
is given by $p_{1}(y|x;\theta)p_{3}(x;\theta)$, where $p_1$ and $p_3$ are given in \eqref{p1} and \eqref{p3}, respectively. Since the $X(C)=x$
is not observed, we integrate it out of the joint density to obtain
\begin{equation}\label{pc2}
P_{C_m}(y;\theta)=P(Y(C)=y,T>C)=\int_{-\infty}^{D}p_{1}(y|x;\theta)p_{3}(x;\theta)dx.
\end{equation}
For a failing item, the joint p.d.f of $Y(T)$ and $T$ equals
\begin{equation}\label{pf}
P_{f}(y,t;\theta)=P(Y(T)=y,T=t<C)=p_{2}(y|t;\theta)f_{T}(t;\theta).
\end{equation}
\subsection{The likelihood}
Assume that $n$ items are on test subject to SSADT over the observation period
$(0,C]$. The sample log-likelihood then is
given by
\begin{eqnarray}
\log L(\theta)&=&\sum_{i=1}^{n}\sum_{j=1}^{m}I_{(\tau_{j-1},\tau_j)}(t_i)\log P_{f}(y_i,t_i;\theta)+(1-I_{(0,C)}(t_i))\log P_{C_m}(y_i;\theta)\nonumber\\
&=&\sum_{j=1}^{m}\sum_{i=\xi_{j-1}+1}^{\xi_j}\log P_{f_j}(y_i,t_i;\theta)
+\sum_{i=\xi_{m-1}+1}^{n}\log P_{C_m}(y_i;\theta).\label{mainloglike}
\end{eqnarray}
$\xi_j=\sum_{k=0}^{j}\nu_k$ and $\nu_0=0$, in which $\nu_j$ is the number of failed
items under stress $S_j$, for $j=1,\ldots,m$,
$(y_{k},t_{k})$, for $k=\xi_{j-1}+1,\ldots,\xi_{j}$, denote the sample failing items
for the stress level $S_j,\; j=1,\ldots,m$, $y_{k}$, for $k=\xi_{m-1}+1,\ldots,n,$ denote the sample surviving (censored) items,
\[P_{fj}(y_i,t_i;\theta)=I_{(\tau_{j-1},\tau_j)}(t_i)\;D(2\pi\sigma_X\sigma_Y)^{-1}(1-\rho^2)^{-1/2}t_i^{-2}e^{-t_i^{-1}Q_{j(y_i,t_i)}},\quad j=1,\ldots,m,\]
\[Q_{j}(y,t)=\eta_1(q_{j}(t,y)-\eta_2P_j(t))^2+\sigma_X^{-2}(D-\mu_{Xj}t)^2/2,,\quad j=1,\ldots,m,\]
\[q_{j}(t,y)=y-\mu_j(y,t),\quad P_{j}(t)=D-\mu_j(x,t),\quad j=1,\ldots,m,\]
in which
\begin{equation}\label{eta12}
\eta_1=\sigma_Y^{-2}(1-\rho^2)^{-1}/2,\quad\mbox{and}\quad
\eta_2=\rho\sigma_Y\sigma_X^{-1}.
\end{equation}
Furthermore, integrating \eqref{pc2} results in
\[P_{C_m}(y_i;\theta)=c_y\sum_{k=1}^{2}(-1)^{k-1}e^{(k-1)\beta_m}\Phi(c_{m}(y;k,1))\phi(c_{m}(y;k,2)),\]
where
\[c_y=\sigma_Y^{-1}C^{-1/2},\quad \beta_m=2D(D-P_m(C))\sigma_X^{-2}C^{-1},\]
\[c_{m}(y;1,1)=\eta_3(P_m(C)-\rho\sigma_X\sigma_Y^{-1} q_{m}(C,y)),\quad c_{m}(y;1,2)=c_yq_{m}(C,y),\]
\[c_{m}(y;2,1)=\eta_3(P_m(C)-\rho\sigma_X\sigma_Y^{-1} q_{m}(C,y)-2D(1-\rho^2)),\quad c_{m}(y;2,2)=c_y(q_{m}(C,y)-2\eta_2D),\]
and
\begin{equation}\label{eta3}
\eta_3=\sigma_X^{-1}(1-\rho^2)^{-1/2}C^{-1/2}.
\end{equation}

\section{Parametric Inference}
In this section, we develop the parametric inferential procedures based on the proposed models.
The maximum likelihood and Bayesian estimation methods are considered for inferential purpose.
From Section 2, it is apparent that the models are analytically intractable. Thus, the finite sample performance
of the maximum likelihood and Bayesian estimators could be examined through a simulation study.
To perform a simulation study,
we set $m=2$, $D=1$, $C=700$,
$S_1=1200$, $S_2=1400$, and $\tau=300,400,500$.
Because of the invariance property of the maximum likelihood estimators, the maximum likelihood estimates of the parameter vector $\theta=(\mu_{X_1},\mu_{X_2},\mu_{Y_1},\mu_{Y_2},\sigma^2_X,\sigma^2_Y,\rho)$ and those of the transformed parameter vector
$$\theta^*=(a,b,c,d,\sigma^2_X,\sigma^2_Y,\rho)$$ can be obtained from each other.
In the following, we assume the transformed parameter vector $\theta^*$ as in Table \ref{tb:margins}.

\begin{table}[!h]
 \footnotesize
\caption{\footnotesize Parameter of model used for the
simulation}\label{tb:margins}
\begin{tabular*}{1\textwidth}%
     {@{\extracolsep{\fill}} |c|c|c|c|c|c|c|c|}
\hline
$\theta^*$&a &b& c & d &$\sigma^2_{X}$&$\sigma^2_{Y}$&$\rho$\\
\hline
  &-2.817991& -4996.008 &-1.644788 &-4995.996&0.001729986&0.0020806801&0.5893698756 \\
\hline
\end{tabular*}
\end{table}
\noindent Using \eqref{mu} we have $(\mu_{X_1},\mu_{X_2},\mu_{Y_1},\mu_{Y_2})=(0.002009813,0.00301472,0.006496424,0.009744636)$.

\subsection{Maximum likelihood}
First, we deal with maximum likelihood estimation of the model parameters.
Suppose $n=30$ independent items are tested subject to SSALT
over the observation period $(0,C]$. The maximum likelihood
estimators (MLEs) of the model parameters can be obtained by
maximizing the log-likelihood \eqref{mainloglike}. It is not possible to obtain the MLEs of
the parameters in a closed form. Thus, numerical computational methods are used for obtaining the MLEs.
A Mont\'{e} Carlo simulation with 10,000 iterations is conducted
using software R 2.14.2 to obtain the estimated relative root of mean square error (RRMSE) and estimated relative
bias (Rbias) of the ML estimators of the parameters. These results are summarized in Table \ref{tb:margins2}.
One can observe from Table \ref{tb:margins2} that the
performance of the estimates are quite satisfactory in terms of RRMSE and Rbias.
\begin{table}[!h]
 \footnotesize
\begin{center}
\caption{\footnotesize Parameter estimates for
$\tau=300,\;400,\;500,\;n=30$}\label{tb:margins2}
\begin{tabular*}{1.0\textwidth}%
     {@{\extracolsep{\fill}} |c|c|ccccccc|}

\hline\hline
$\tau$ & & $\mu_{X_1}$ & $\mu_{X_2}$& $\mu_{Y_1}$ & $\mu_{Y_2}$ & $\sigma^2_{X}$ & $\sigma^2_{Y}$&$\rho$ \\
\hline
 &MLE     &0.001544& 0.002315 &0.006268 &0.009401&0.001755 &0.002020&0.591156 \\
300&Rbias&-0.231723&-0.154533&-0.035092&-0.035226&0.014603&-0.029246&0.003031  \\
 &RRMSE   &0.263621 &0.175785 &0.063134 &0.063224&0.271781 &0.245912&0.204330\\
 \hline
 \hline
 &MLE     &0.001536& 0.002303 &0.006277 &0.009415&0.001839 &0.002045&0.596587\\
400&Rbias&-0.236005&-0.157394&-0.033721&-0.033851&0.062744&-0.017047&0.012246 \\
 &RRMSE   &0.275125 &0.183462 &0.063961 &0.064030&0.271558 &0.248880&0.203052  \\
 \hline
 \hline
 &MLE    & 0.001599&0.002396 &0.006324 &0.009484 &0.001850 &0.002050&0.596866\\
500&Rbias&-0.204660&-0.13650&-0.026547&-0.026784 &0.069158&-0.014685&0.012720  \\
 &RRMSE  & 0.262404&0.174980 &0.063265 &0.064235 &0.273821 &0.246177&0.204354  \\
 \hline
\end{tabular*}
\end{center}
\end{table}

\begin{table}[!h]
\footnotesize
\begin{center}
\caption{\footnotesize Thirty simulated observation for parameter
set of Table \ref{tb:margins}. }\label{thirty}
\begin{tabular*}{0.75\textwidth}%
     {@{\extracolsep{\fill}}|ccc|ccc|ccc|}
     \hline\hline
& $\tau=300$ &  &&$\tau=400$&  &&$\tau=500$&  \\
\hline
 \hline
   $\delta$& t&y&$\delta$&t &y&$\delta$&t&y \\
\hline
1 &206 &2.9043836&1&125&1.3302026 &1&72&0.8630739 \\
1 &204&2.2834415 &1&347&2.0791452&3&700&5.7413420 \\
2&358 &2.0369846 &2&409&2.8987105&1&257&0.5140461\\
2&424&2.2286551 &3&700&4.3707818&2&627&3.7818455 \\
2&528&3.3882536&1&321&4.2948213&1&265&2.5664485\\
1&293&1.0765821&2&664&5.1428573&3&700&2.3397970\\
2&433&3.4253562&2&413&2.1084361&2&588&4.4356021\\
2&367&2.7105020&2&575&3.8347019&1&261&2.2925127\\
2&481&2.7411018&3&700&4.3353895&1&152&1.5052757\\
1&74&0.4584009&1&61&1.4195226&1&203&2.2968604\\
1&232&1.4229018&2&443&4.7402742&3&700&5.2271200\\
2&563&2.0737839&1&74&0.9594538&1&205&2.4270830\\
2&524&4.6559941&3&700&6.6440064&2&500&3.9215091\\
2&398&3.0469754&2&439&2.2434726&3&700&4.2650212\\
1&83&1.1645206&2&543&3.6592403&2&521&2.0058003\\
1&288&1.5370298&3&700&3.9542721&1&321&3.1932579\\
2&518&2.8000903&1&238&0.4759539&1&435&3.1052309\\
2&558&4.4736314&1&104&1.2579545&1&160&2.6871790\\
1&106&1.3271670&2&413&2.6969144&1&329&3.0110215\\
1&699&7.4817986&2&429&1.5348759&2&687&5.7361510\\
2&538&4.5781005&1&231&0.9987282&1&249&1.9830004\\
1&98&0.3647197&1&205&2.1099217&2&578&2.8660246\\
1&184&1.5738009&3&700&6.1866970&1&335&2.4279700\\
2&379&2.8413248&1&146&1.8776734&1&273&2.3275245\\
1&102&1.1580797&3&700&4.6567010&1&143&1.3927807\\
1&165&1.6696197&1&375&1.5709192&1&161&2.3288032\\
2&584&3.6045384&2&541&2.8358232&2&692&5.5437838\\
2&371&2.4435304&2&600&4.1900476&1&175&0.7888585\\
2&538&4.4705936&2&623&3.8115301&1&199&1.0319427\\
2&303&2.8150623&1&274&2.0937201&3&700&3.1902784\\
 \hline
   \end{tabular*}
  \end{center}
    \end{table}
\subsection{Bayesian approach}
The Bayesian approach is
appealing to statisticians and reliability engineers,
since it provides a method of using their past experiences
and/or prior convictions for inference. From a Bayesian point of
view, we can treat the unknown parameters as a random variable with
a known prior probability distribution. Then, we can combine
information from the random sample and prior probability
distribution to obtain the Bayesian estimators for the parameters of
the model. However, in most practical applications, where the Bayesian
approach is used, it is difficult to compute analytically the
posterior distribution. The Markov chain Mont\'{e} Carlo (MCMC) method uses to generate a
sample from the posterior distribution large enough so that any
desired feature of the posterior distribution can be accurately
obtained. Because of the restrictions $\mu_{X_1}<\mu_{X_2}$ and
$\mu_{Y_1}<\mu_{Y_2}$, we have to consider joint priors for the vectors
$(\mu_{X_1},\mu_{X_2})$ and $(\mu_{Y_1},\mu_{Y_2})$, while we can
consider independent priors for the transformed parameters $a,b,c$ and $d$.
To simplify the calculations, we perform the Bayesian approach for the transformed
parameter vector $\theta^*=(a,b,c,d,\sigma^2_X,\sigma^2_Y,\rho)$.

\begin{table}[!h]
\footnotesize
\begin{center}
\caption{\footnotesize Parameter estimation results for
$\tau=300,\;400,\;500,\;n=30$}\label{std}
\begin{tabular*}{0.93\textwidth}%
     {@{\extracolsep{\fill}}|cc|cccccc|}
 \hline\hline
 & &  &   &$\tau=300$ &  &  &  \\
\hline
 && Mean & Std& MC-er &2.5$\%$ &Medain &97.5$\%$ \\
\hline
 &a&-3.21803700&0.1809108&0.0009045538&-3.562406&-3.185069&-2.991624 \\
&b&-4100.17100&317.9782&1.589891&-4624.464&-4109.659&-3580.026 \\
 &c&-2.24168800&0.2151238&0.001075619&-2.586862&-2.256567&-1.897493 \\
 &d&-4085.98400&353.0725&1.765363&-4670.439&-4076.054&-3539.374 \\
&$\sigma^2_X$&0.001651598&0.00001052115&5.260574$\times10^{-8}$&0.001635779&0.001649452&0.001679365  \\
 &$\sigma^2_Y$ &0.002155083&0.00007563872&3.781936$\times10^{-7}$&0.00203197&0.002152928&0.002293804  \\
  &$\rho$&0.594830600&0.03566883&1.783442$\times10^{-4}$&0.5374839&0.5871186&0.6582097 \\
  \hline\hline
& &   & & $\tau=400$ &  &  &  \\
\hline
 && Mean & Std& MC-er &2.5$\%$ &Medain &97.5$\%$ \\
\hline
 &a&-3.318978&0.1811505&0.0009057527&-3.599276&-3.355476&-3.028328 \\
&b&-4070.080&300.5802&1.502901&-4603.892&-4046.055&-3633.086 \\
 &c&-2.245963&0.1760872&0.0008804360&-2.539122&-2.250249&-1.975724 \\
 &d&-4102.815&301.8582&1.509291&-4586.286&-4105.338&-3627.138 \\
&$\sigma^2_X$&0.001817105&5.632395$\times10^{-5}$&2.816198$\times10^{-7}$&0.001746018&0.001798002&0.001932355  \\
 &$\sigma^2_Y$ &0.001924476&3.424275$\times10^{-5}$&1.712138$\times10^{-7}$&0.001872394&0.001920889&0.001980292  \\
  &$\rho$&0.5896328&0.01278770&6.393848$\times10^{-5}$&0.5574416&0.5901468&0.6118633 \\
\hline\hline
&&   &  & $\tau=500$ &  &  &  \\
\hline
& & Mean & Std& MC-er &2.5$\%$ &Medain &97.5$\%$ \\
\hline
 &a&-3.204618&0.1443177&0.0007215886&-3.447255&-3.158764&-2.968463 \\
&b&-4151.146&312.5983&1.562992&-4692.574&-4130.741&-3667.357 \\
 &c&-2.250095&0.1716518&0.0008582590&-2.545561&-2.212597&-1.971093 \\
 &d&-4090.401&288.7204&1.443602&-4581.323&-4123.713&-3607.695 \\
&$\sigma^2_X$&0.001729962&0.00005117870&2.558935$\times10^{-7}$&0.001638141&0.001740245&0.001804733  \\
 &$\sigma^2_Y$ &0.002047766&0.00003888018&1.944009$\times10^{-7}$&0.001976453&0.002049934&0.002135915  \\
  &$\rho$&0.5792594&0.01128655&5.643275$\times10^{-5}$&0.5586604&0.5785074&0.6028768 \\
  \hline
  \end{tabular*}
\end{center}
\end{table}

\begin{table}[!h]
 \footnotesize
\begin{center}
\caption{\footnotesize Parameter estimates for non informative prior
$\tau=300,\;400,\;500,\;n=30$}\label{rbias}
\begin{tabular*}{1.0\textwidth}%
     {@{\extracolsep{\fill}} |c|c|ccccccc|}
\hline\hline
& &  &  & & $n=30$ &  &  &  \\
\hline
$\tau$ & & $\mu_{X_1}$ & $\mu_{X_2}$& $\mu_{Y_1}$ & $\mu_{Y_2}$ & $\sigma^2_{X}$ & $\sigma^2_{Y}$&$\rho$ \\
\hline
 &BE      &0.002479&0.003456&0.006636 &0.009242 &0.001652&0.002155&0.594831\\
300&Rbias &0.233561&0.146273&0.021469 &-0.05154&-0.045311&0.035759&0.009265 \\
 &MC-er&7.27$\times10^{-7}$&8.05$\times10^{-7}$&8.73$\times10^{-7}$&1.31$\times10^{-5}$&1.05$\times10^{-5}$&7.56$\times
10^{-5}$&0.000178 \\
 \hline
 \hline
 &BE     &0.002285&0.003178&0.006534 &0.009112&0.001817 &0.001924&0.589633\\
400&Rbias&0.136740&0.054138&0.005756&-0.064964&0.050358&-0.075074&0.000446 \\
 &MC-er&3.80$\times10^{-7}$&3.61$\times10^{-7}$&1.08$\times10^{-6}$&1.29$\times10^{-5}$&2.82$\times10^{-7}$&
 1.71$\times10^{-6}$&6.39$\times10^{-5}$  \\
 \hline
 \hline
 &BE      &0.002430&0.003399&0.006568 &0.009141 &0.001730 &0.002048  &0.579259\\
500&Rbias &0.209242&0.127457&0.009893&-0.006195&-0.000014&-0.015819 &-0.017155 \\
 &MC-er&9.35$\times10^{-7}$&9.50$\times10^{-7}$&8.70$\times10^{-7}$&1.29$\times10^{-5}$&2.56$\times10^{-7}$&1.94$\times10
 ^{-7}$&5.64$\times10^{-5}$  \\
 \hline
\end{tabular*}
\end{center}
\end{table}

Table \ref{thirty} presents simulated data sets by using the parameters in
Table \ref{tb:margins} for $\tau=300,400,500$. We consider the Bayes estimation
of the transformed parameter vector, ${\theta}^*$, based on data sets in Table \ref{thirty},
under the square error and absolute error loss
functions. An analytic calculation of estimators and their risks for comparison
is far from reach. To carry out an empirical comparison, a simulation study was conducted
using software R 2.14.2 to generate a sequence of parameter values from the posterior density
of $\theta^*$ given the generated data set of Table \ref{thirty} by making use of the random walk
Metropolis-Hasting algorithm.

\begin{figure}[!hbtp]
\centerline{\psfig{figure=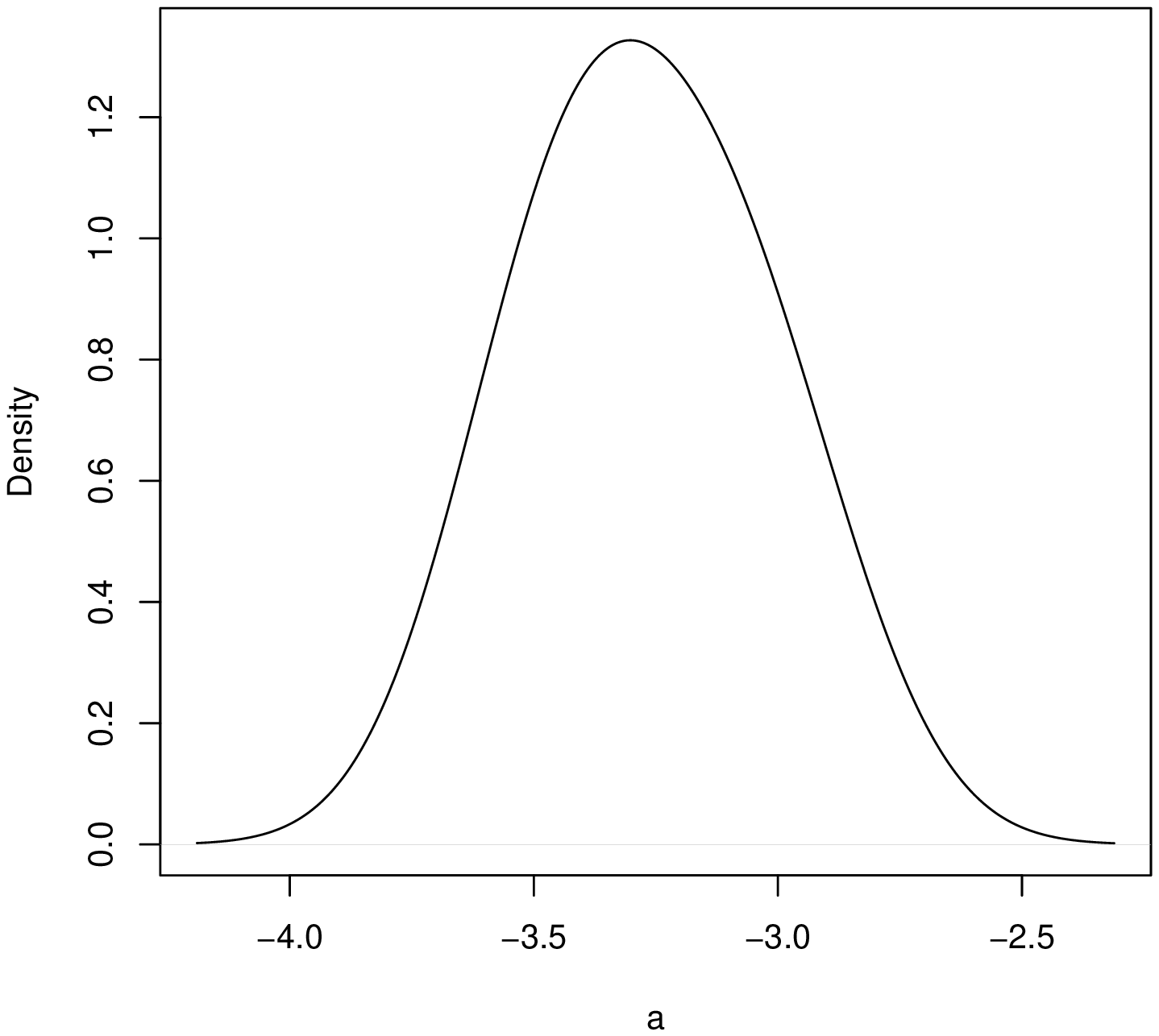,width=4cm}\psfig{figure=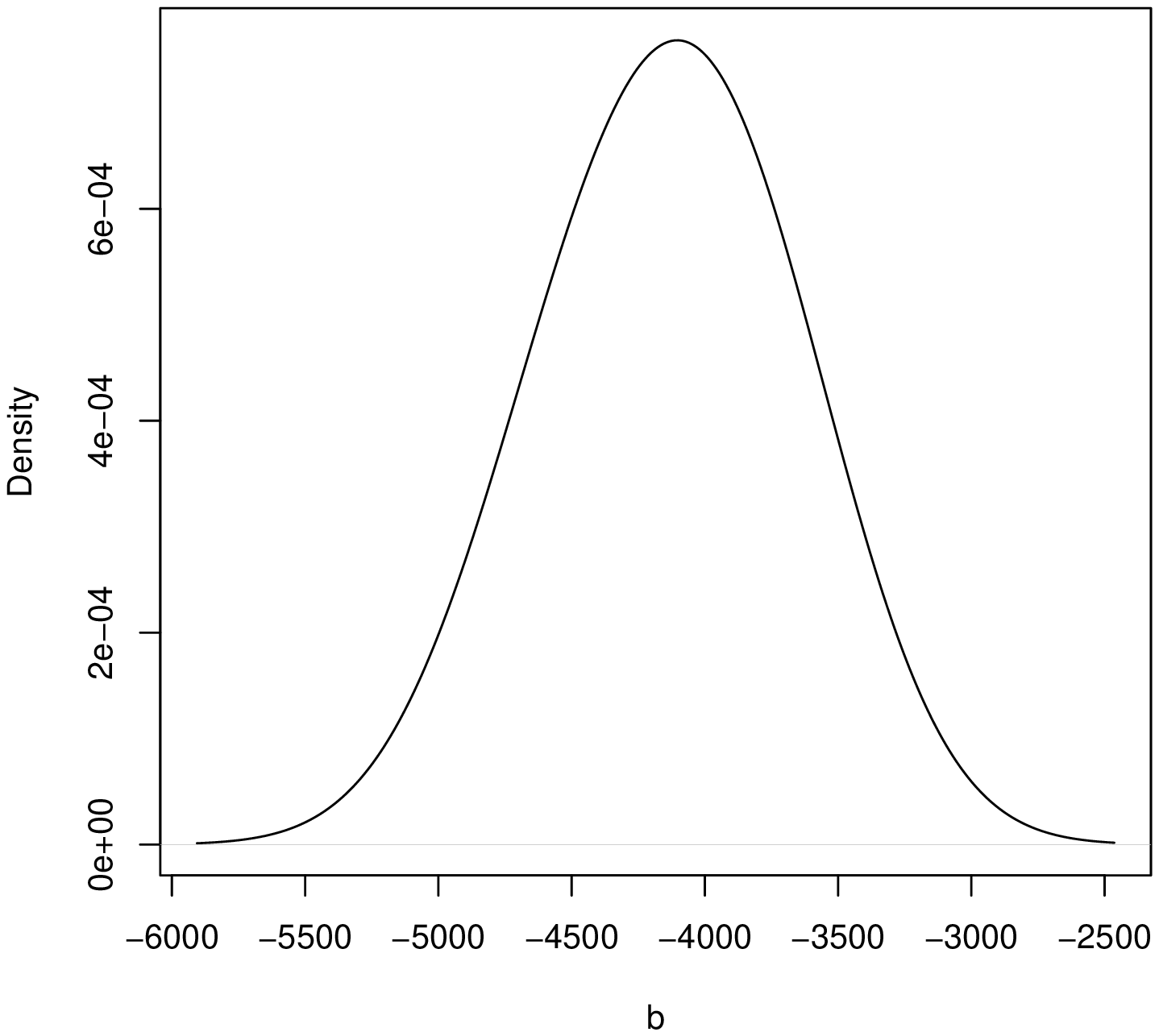,width=4cm}\psfig{figure=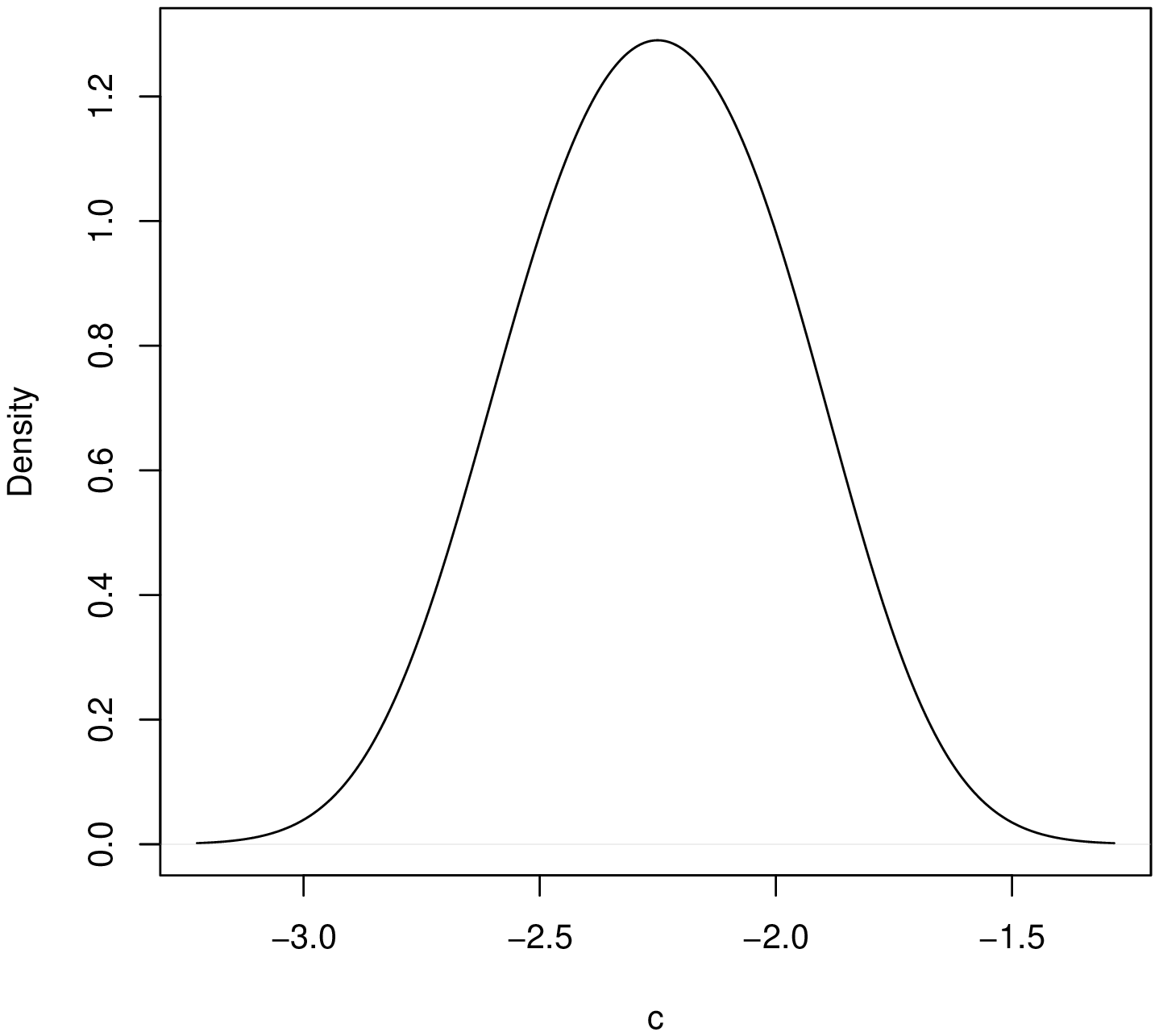,width=4cm}}
\centerline{\psfig{figure=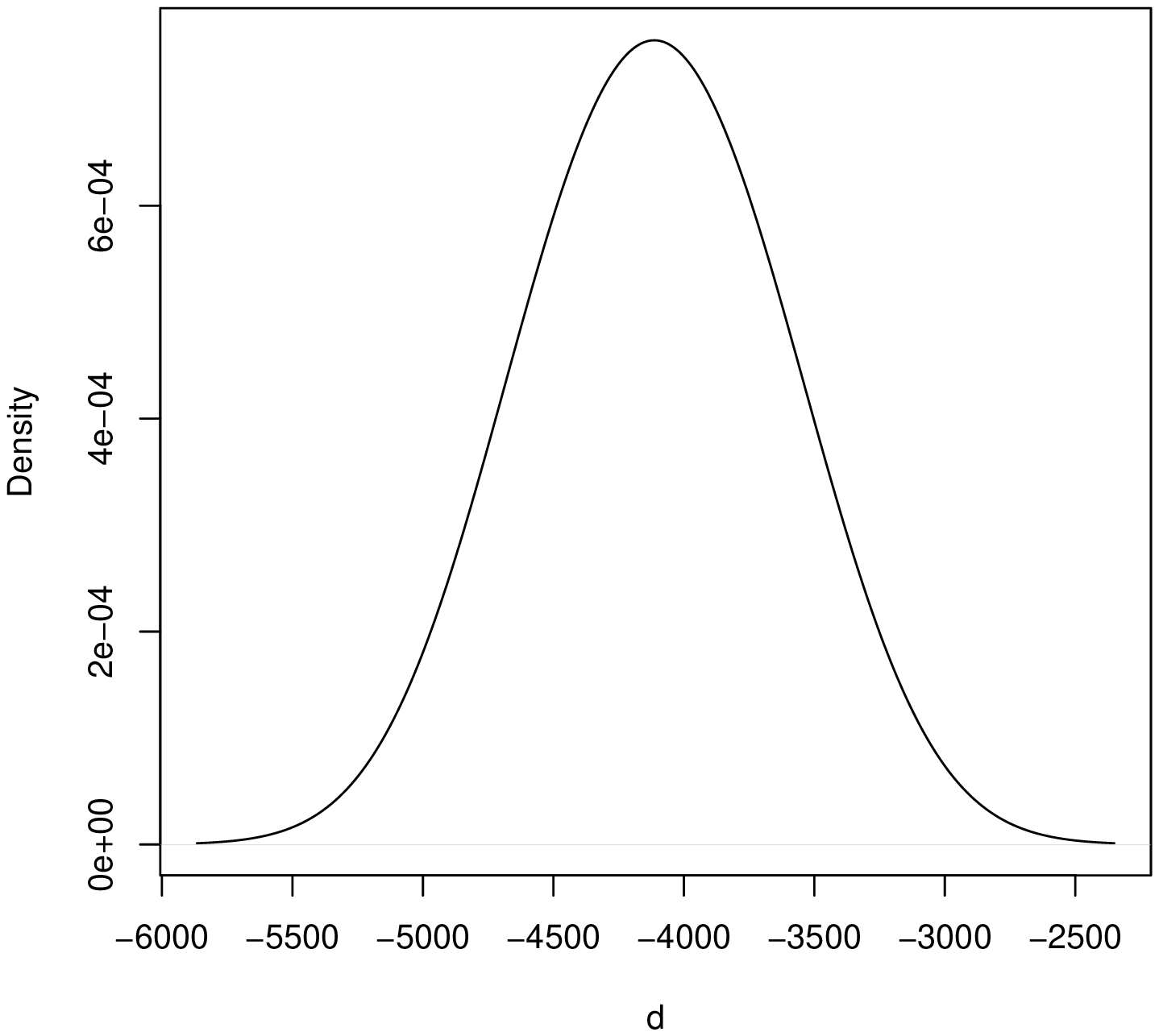,width=3cm}\psfig{figure=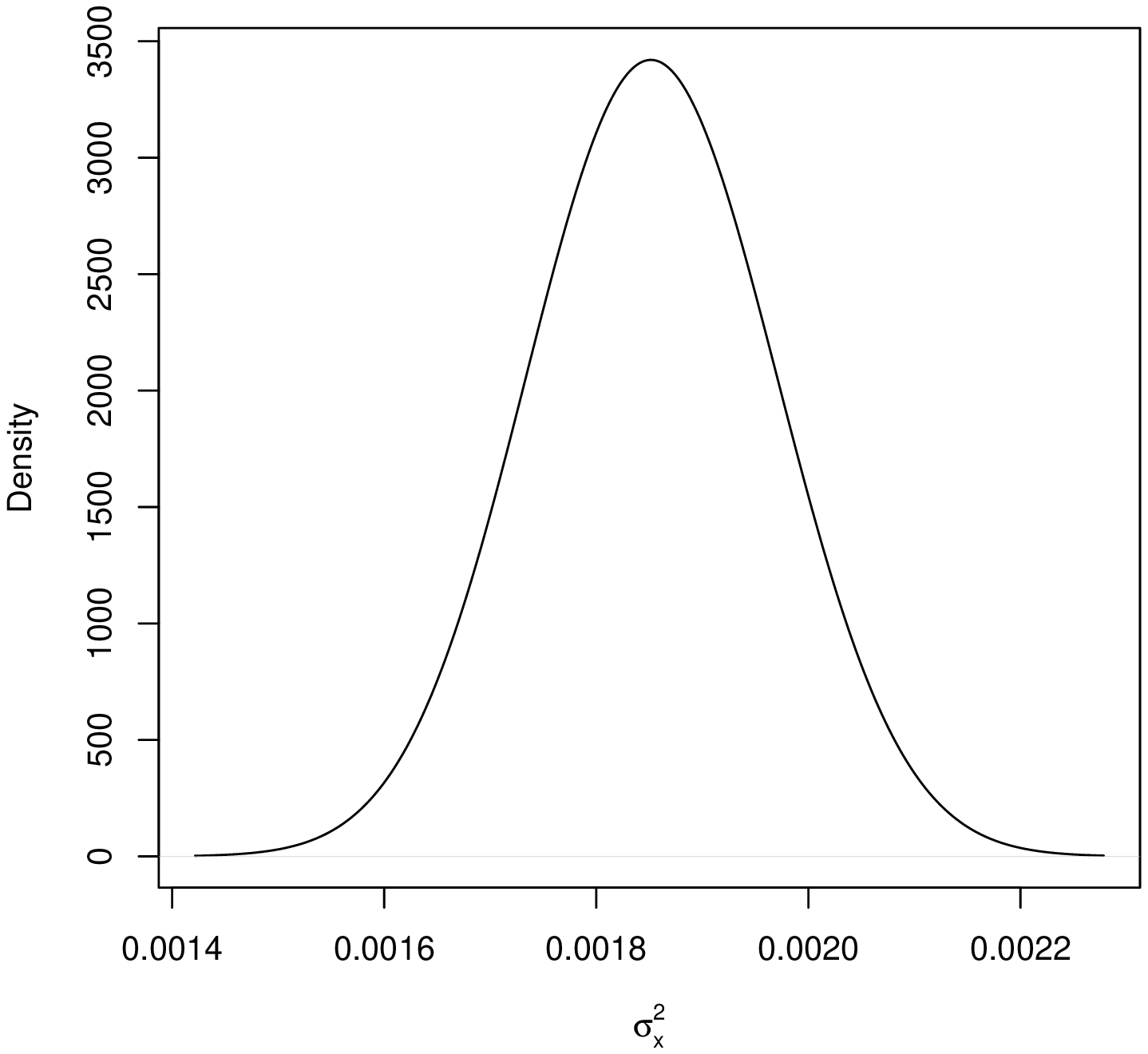,width=3cm}\psfig{figure=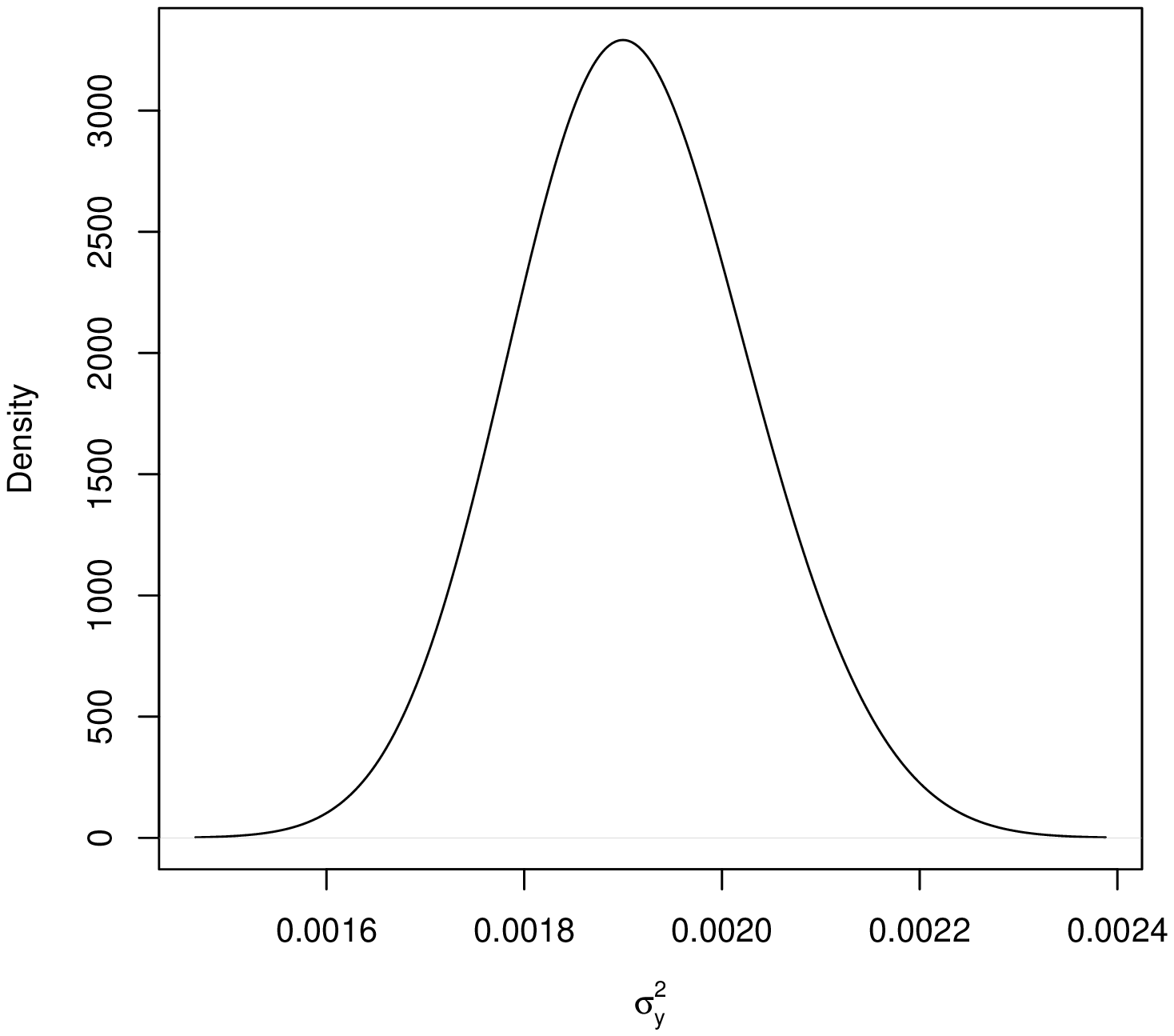,width=3cm}\psfig{figure=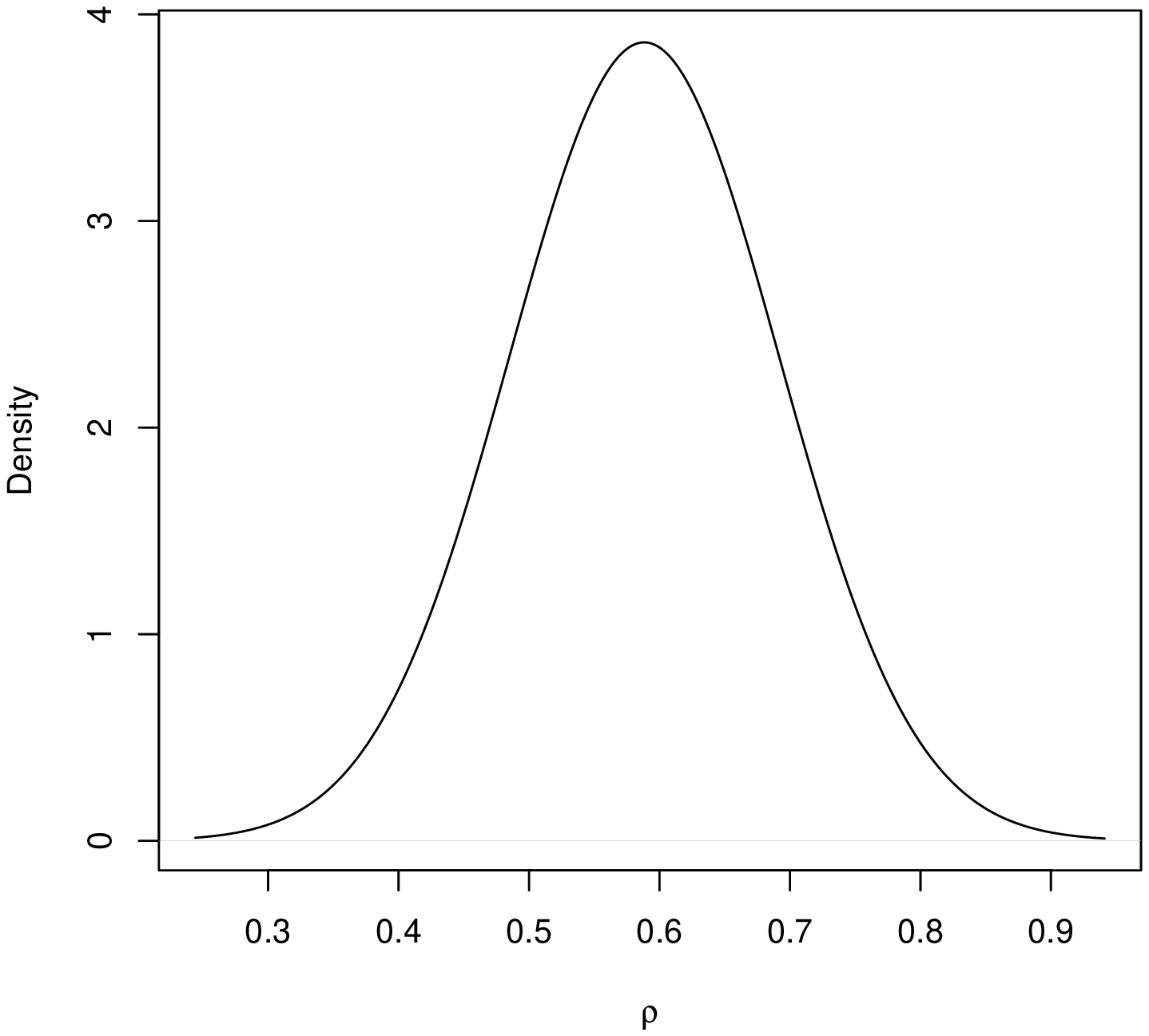,width=3cm}}
 \caption{\scriptsize The empirical posterior densities of the model parameters for $\tau=400$.}\label{f1}
\end{figure}

\begin{figure}[!hbtp]
\centerline{\psfig{figure=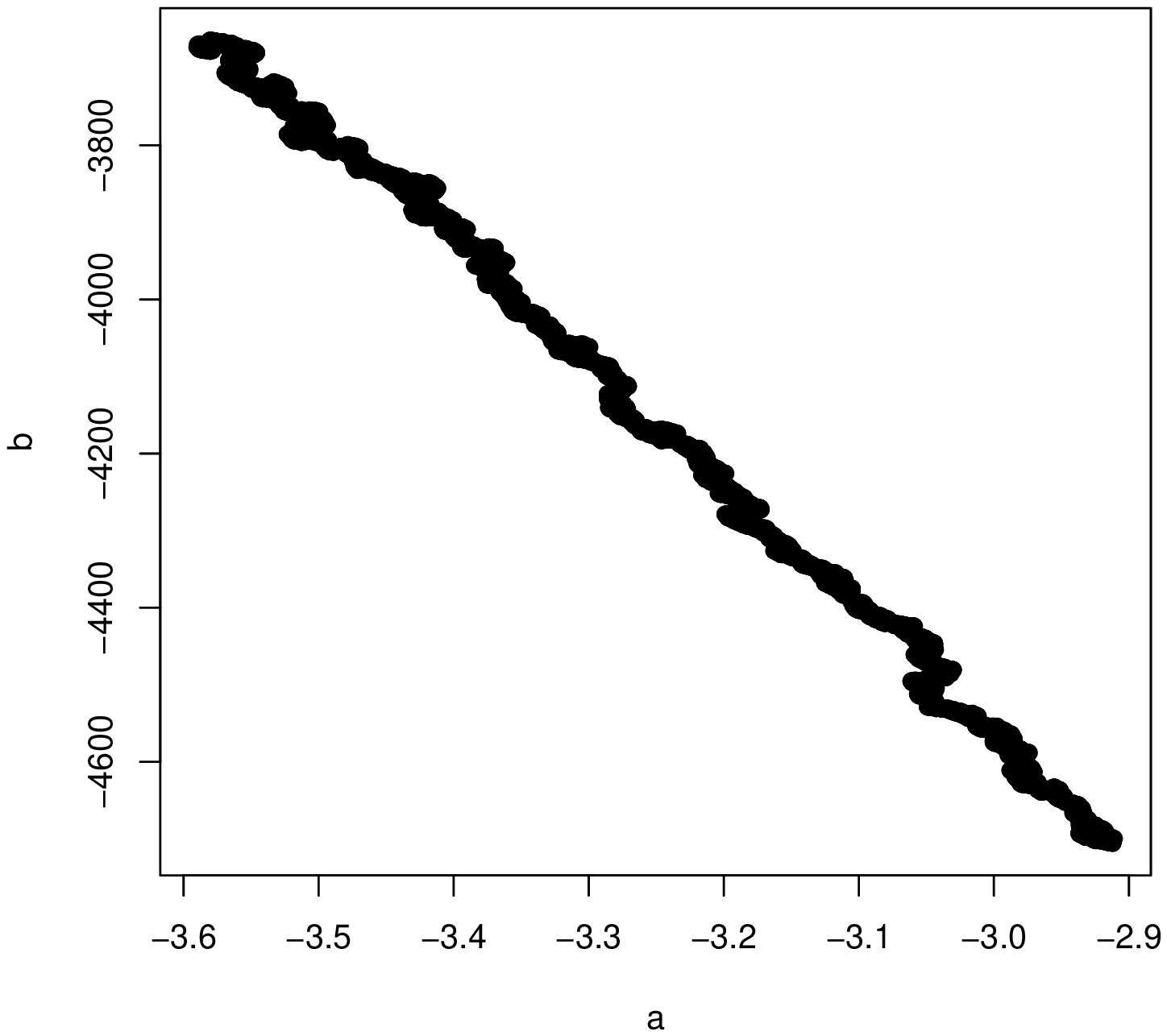,width=4cm}\psfig{figure=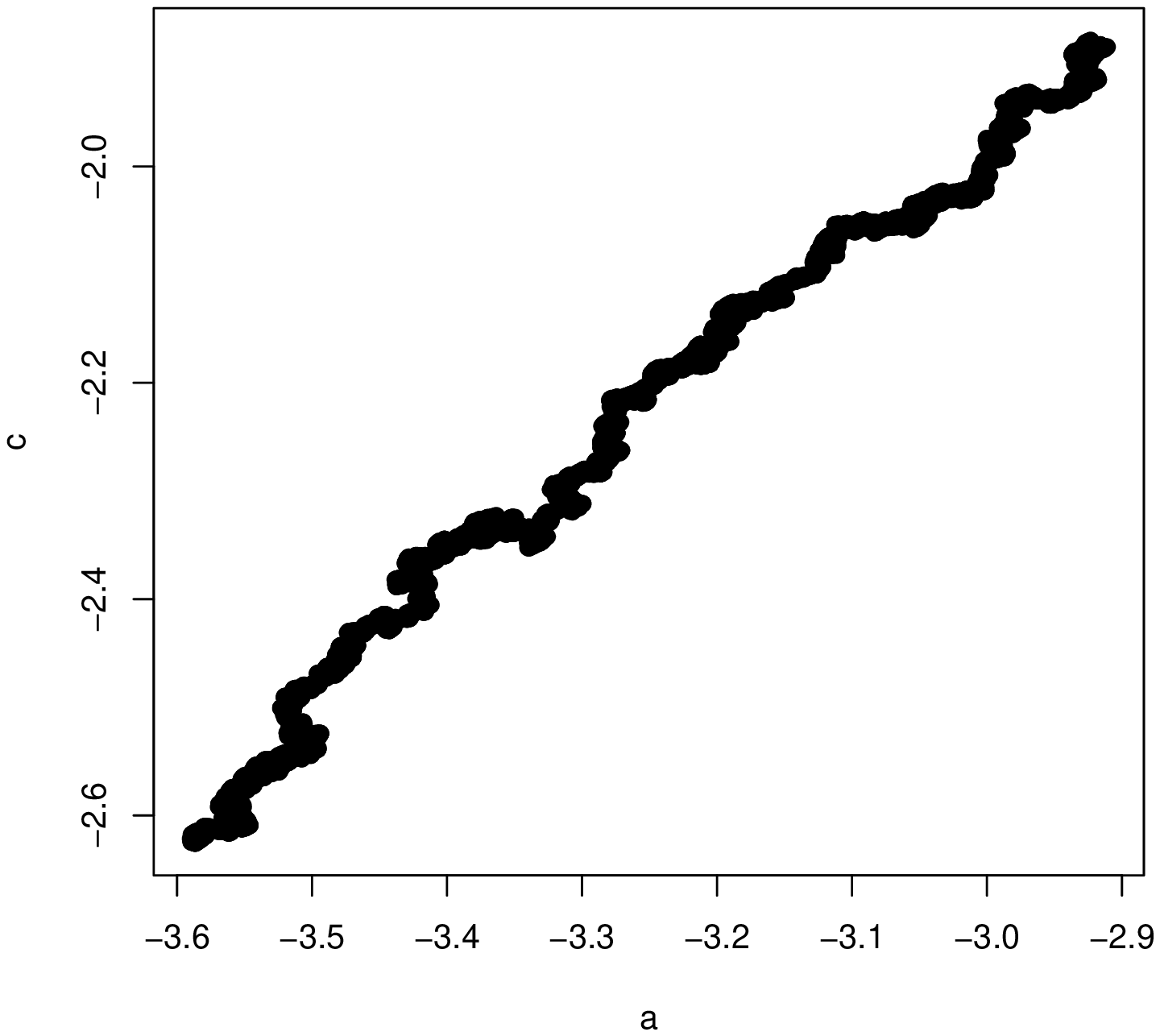,width=4cm}\psfig{figure=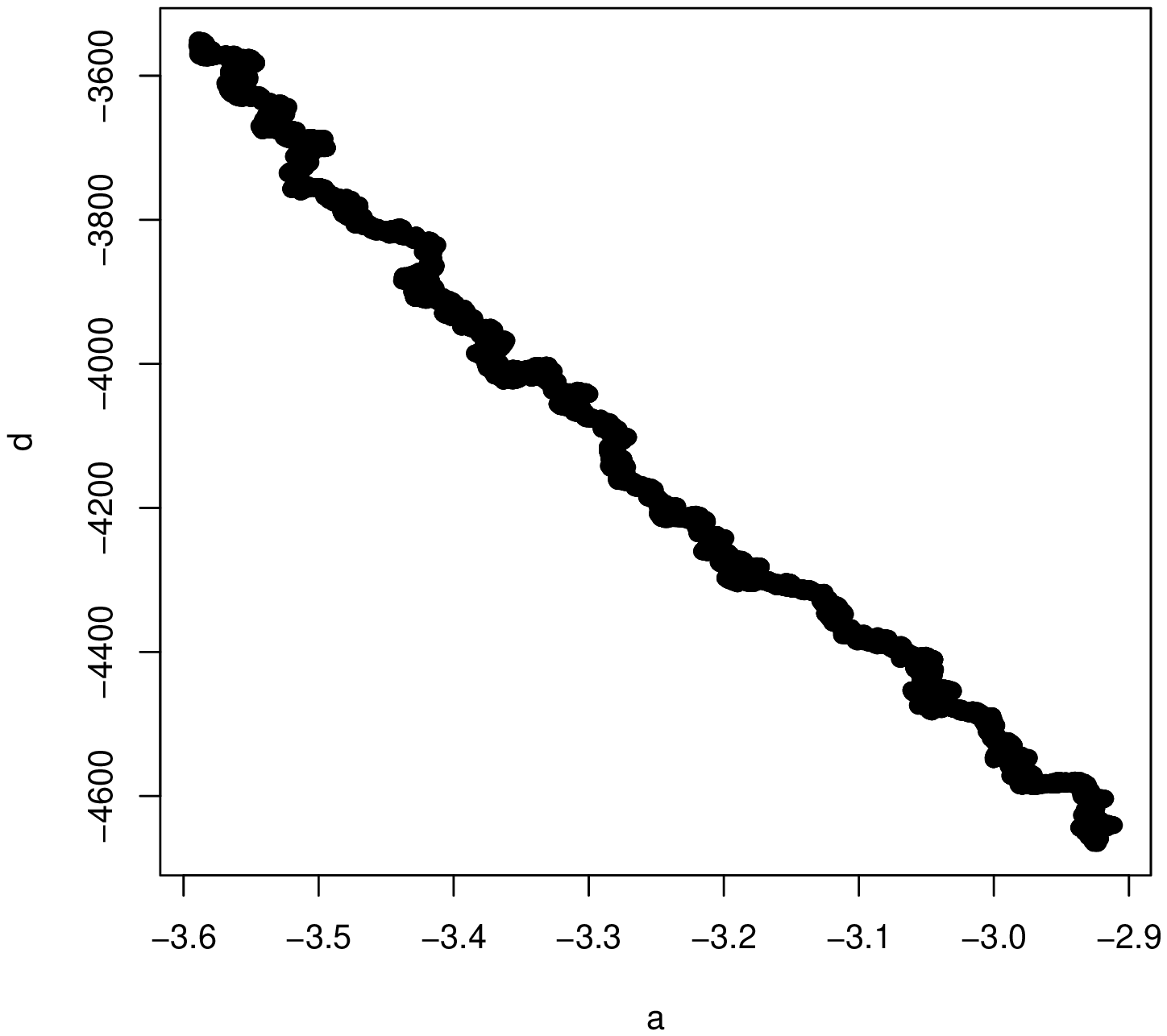,width=4cm}}
\centerline{\psfig{figure=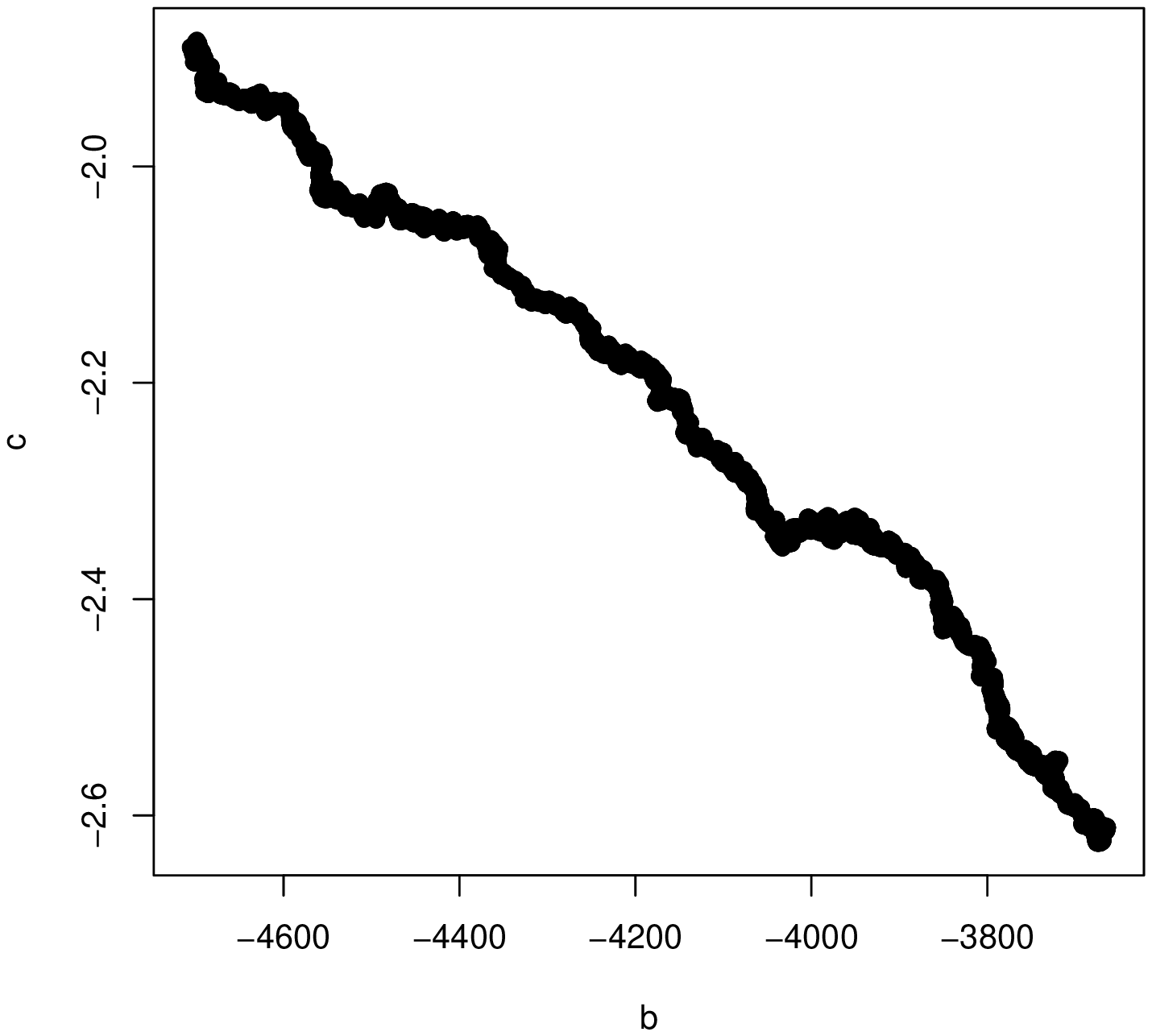,width=4cm}\psfig{figure=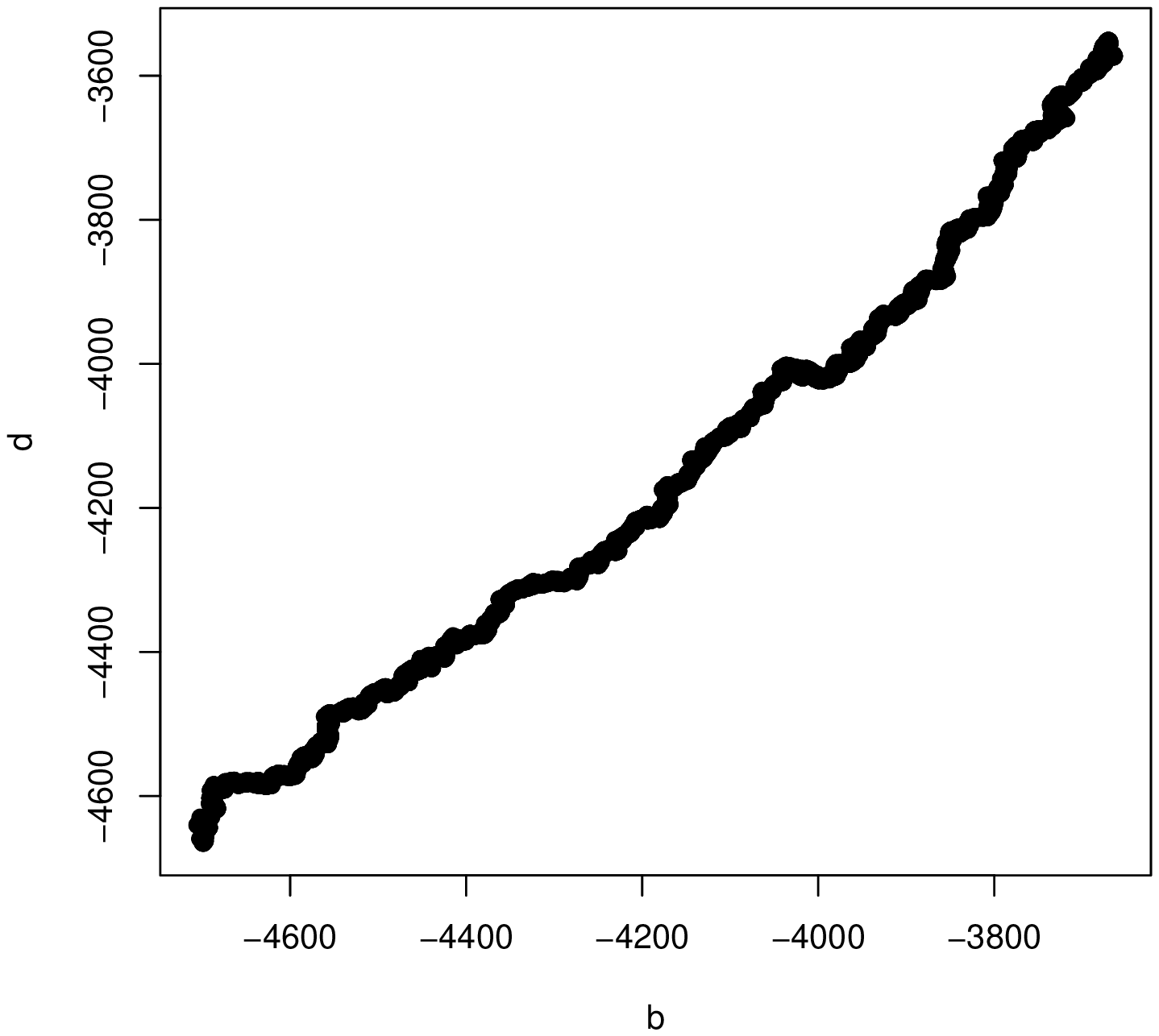,width=4cm}\psfig{figure=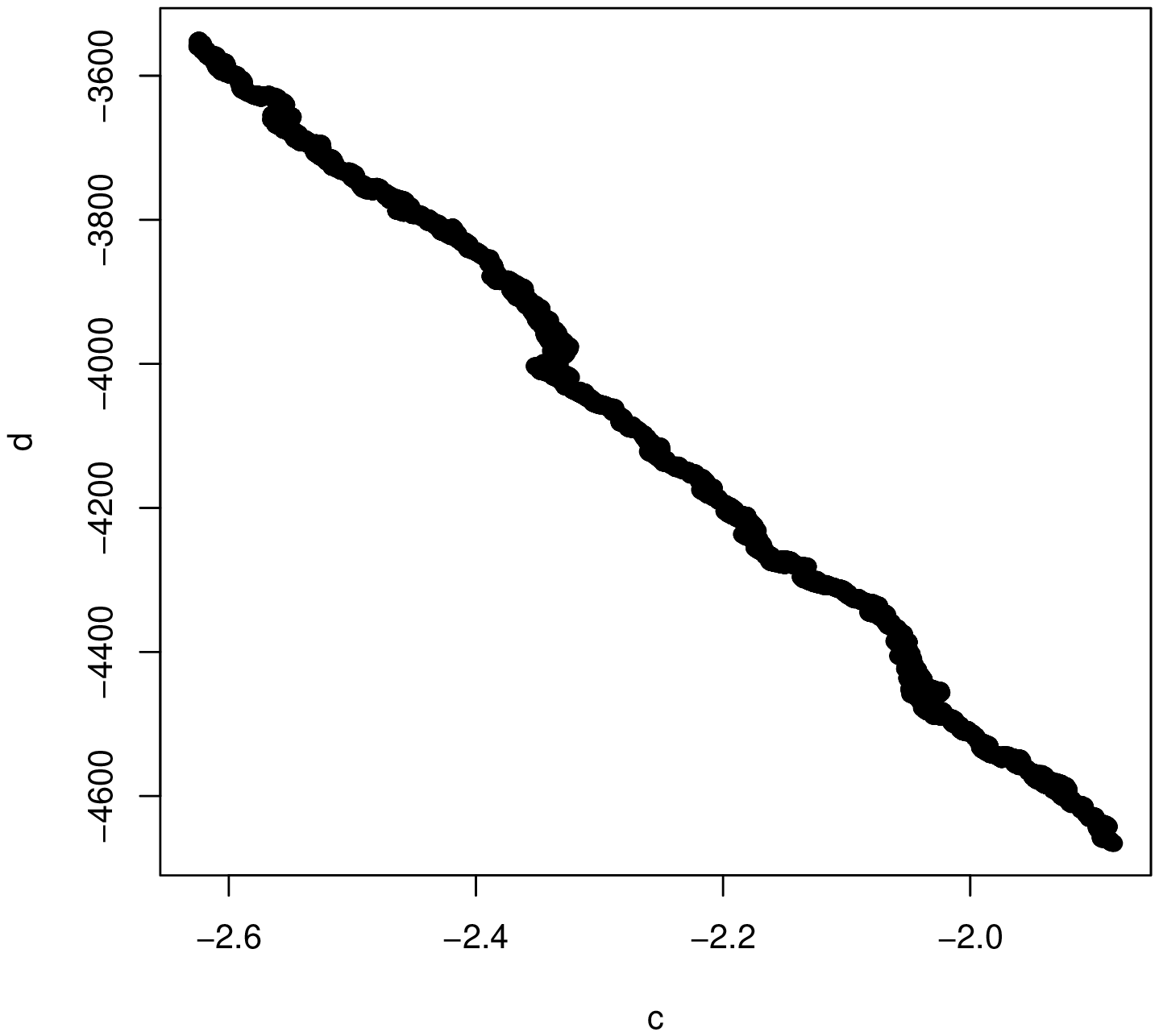,width=4cm}}
 \caption{\scriptsize Two dimensional plots of the generated model parameters $(a,b,c,d)$ for $\tau=400$.}\label{f2}
\end{figure}

To facilitate the Bayesian approach, we assume independent prior distributions for the model
parameters, that is
$$\pi(a,b,c,d,\sigma^2_{X},\sigma^2_{Y},\rho)\propto \pi_1(a)\pi_2(b)\pi_3(c)\pi_4(d)
 \pi_5(\sigma^2_{X})\pi_6(\sigma^2_{Y})\pi_7(\rho),$$
where $\pi_1(a),\;\pi_2(b),\;\pi_3(c)$ and $\pi_4(d)$ are assumed to be the low informative normal densities with zero mean and the variance equal to $10^4$,
$\pi_5(\sigma^2_{X})$ and $\pi_6(\sigma^2_{Y})$ are assumed to be the non-informative Jeffrey's priors $\pi(\sigma^2_{X})\propto\frac{1}{\sigma^2_{X}},\;\pi(\sigma^2_{Y})\propto\frac{1}{\sigma^2_{Y}}$ and $\pi_7(\rho)$ is taken to be the non-informative
uniform$(-1,1)$ prior.

The random walk Metropolis-Hasting algorithm is executed 50000 times and the last 40000 were used for the sake of convergency.
The empirical posterior densities of the model parameters and two dimensional plots of the generated model parameters $(a,b,c,d)$
are shown in Figures \ref{f1} and \ref{f2}, respectively, for $\tau=400$.  Using these empirical densities we estimate the
mean, standard deviation (Std), MCMC error (MC-er), the median and other critical quantiles of parameters.
These numerical results are summarized in Table \ref{std}.

For the sake of brevity, only the values of the Bayse estimates (BEs)
based on the square error loss, as well as their Rbias, MC-er were typically given in Table
\ref{rbias} to be compared with the corresponding values of the ordinary MLEs. Similar comparisons can be made between
BEs based on the absolute error loss and the MLEs.

\section{Optimal test plan}

For $m=2$ stress levels, we have
\begin{equation}\label{linear}
 \mu_{X0}=\exp\left(\frac{\log\mu_{X1}-\alpha\log\mu_{X2}}{1-\alpha}\right)\\
 \mbox{and}\\
 \mu_{Y0}=\exp\left(\frac{\log\mu_{Y1}-\alpha\log\mu_{Y2}}{1-\alpha}\right),\\
 \end{equation}
 where
 $$\alpha=\frac{S_1-S_0(273+S_2)}{S_2-S_0(273+S_1)}$$
 is called the {\it stress ratio}.

The optimization criterion considered in this paper has to find the
optimal stress changing time $0<\tau^*<C$ which minimizes the Approximate variance (Avar)
of the ML estimate of the 100$p^{\mbox{th}}$ percentile of the distribution of $T$, $\hat{\xi}_p$, under the
normal stress level $S_0$. The Avar of $\hat{\xi}_p$ is a function of the stress changing time
$\tau$ and the parameter vector $\theta$. Hence, before performing the
optimization procedure, one have to estimate the parameter vector
$\theta$ using a lifetime data in normal conditions. This is done
via the ML estimation using \eqref{mainloglike} and based on a pilot study.

The Avar of $\hat{\xi}_p$ can be obtained
as a function of the approximated variance of the MLE of $\theta$ ( the inverse of the Fisher information matrix,
 $I(\theta)$), using the delta method as
$$\mbox{Avar}(\hat{\xi}_p)=H'I^{-1}(\theta)H/(f_{T_0}(\hat{\xi}_p))^2,$$
where $f_{T_0}(t)$ is the corresponding pdf of $G_0(t)$ in \eqref{gt} and
\[H'=\left[\frac{\partial\hat{{G}}_0(\hat{\xi}_p)}{\partial\hat{\mu}_{X1}},
\frac{\partial\hat{{G}}_0(\hat{\xi}_p)}{\partial\hat{\mu}_{X2}},0,0,
\frac{\partial\hat{{G}}_0(\hat{\xi}_p)}{\partial\hat{\sigma}_X^2},0,0\right].\]
Note that ${{G}}_0(s)$ is not a function of
${\mu}_{Y1},{\mu}_{Y2},{\sigma}_Y^2$, and $\rho$.

We have
\begin{eqnarray*}
\frac{\partial\hat{{G}}_0(\hat{\xi}_p)}{\partial \hat{\mu}_{X_1}
}&=&\frac{\hat{\xi}_p\hat{\mu}_{X_0}\phi\left(c_{1x}\right)}{\hat{\mu}_{X_1}(1-\alpha)\sqrt{\hat{\sigma}_X^2\hat{\xi}_p}}
+\frac{2D\hat{\mu}_{X_0}e^{\beta_3}}{\hat{\mu}_{X_1}\hat{\sigma}_X^2(1-\alpha)}\Phi\left\{c_{2x}\right\}-
\frac{\hat{\xi}_p\hat{\mu}_{X_0} e^{\beta_3}}{\hat{\mu}_{X_1}(1-\alpha)\sqrt{\hat{\sigma}_X^2\hat{\xi}_p}}\phi\left\{c_{2x}\right\},\\
\end{eqnarray*}
\[\frac{\partial\hat{{G}}_0(\hat{\xi}_p)}{\partial \hat{\mu}_{X_2}}=
-\alpha\frac{\hat{\mu}_{X_1}}{\hat{\mu}_{X_2}}\frac{\partial\hat{{G}}_0(\hat{\xi}_p)}{\partial \hat{\mu}_{X_1}},\]
and
\begin{eqnarray*}
\frac{\partial\hat{{G}}_0(\hat{\xi}_p)}{\partial \hat{\sigma}_{X}^2}
&=&-\frac{\left[\hat{\mu}_{X_0}\hat{\xi}_p-D\right]\phi\left(c_{1x}\right)}{2\hat{\sigma}_X^3\sqrt{\hat{\xi}_p}}
-\frac{2D\hat{\mu}_{X_0}e^{\beta_3}\Phi\left\{c_{2x}\right\}}{\hat{\sigma}_X^4}+\frac{\left[\hat{\mu}_{X_0}\hat{\xi}_p+D\right]e^{\beta_3}\phi
\left\{c_{2x}\right\}}{2\hat{\sigma}_X^3\sqrt{\hat{\xi}_p}}.
\end{eqnarray*}
where
\[c_{1x}=\sqrt{\frac{1}{\hat{\sigma}_X^2\hat{\xi}_p}}\left[\hat{\mu}_{X_0}\hat{\xi}_p-D\right],\quad c_{2x}=-\sqrt{\frac{1}{\hat{\sigma}_X^2\hat{\xi}_p}}\left[\hat{\mu}_{X_0}\hat{\xi}_p+D\right],\]
\[\beta_3=\frac{2D}{\hat{\sigma}_X^2}\hat{\mu}_{X_0},\quad \hat{\mu}_{X_0}=\exp\left(\frac{\log\hat{\mu}_{X_1}-\alpha\log\hat{\mu}_{X_2}}{1-\alpha}\right)\]
and $\hat{\mu}_{X1}$, $\hat{\mu}_{X2}$ and $\hat{\sigma}_X^2$ are
the MLEs of $\mu_{X1}$, $\mu_{X2}$ and $\sigma_X^2$, respectively,
which are computed numerically using the log-likelihood in \eqref{mainloglike}.

In order to calculate the estimate of the Fisher information matrix
of the data at $\theta$, that is $I(\theta)=((I_{r,s}(\theta)))$, first let
$\theta=(\mu_{X1},\mu_{X2},\mu_{Y1},\mu_{Y1},\sigma_X^2,\sigma_Y^2,\rho)=(\theta_1,\theta_2,\ldots,\theta_7)$.
The random vector ($\nu_1,\nu_2$) in \eqref{mainloglike} follows a multi-nomial distribution
with parameters $n,\;p_1=G_1(\tau_1),\;p_2=G_2(C)-G_2(\tau_1)$, where
$G_j(t),\;j=1,2$ are as in \eqref{gt} with $\mu_{X_0}$ replaced by $\mu_{X_j},\;j=1,2$ respectively.
We have
\begin{eqnarray*}
\hat{I}_{r,s}(\hat{\theta})&=&E\left.\left(E\left(\left.\frac{-\partial^2\log
L(\theta)}{\partial\theta_r\partial\theta_s}\right|\nu_1,\nu_2\right)\right)\right|_{\hat{\theta}}\\
&=&\sum_{\nu_1=0}^n
\sum_{\nu_2=0}^{n-\nu_1}{n \choose \nu_1}{n-\nu_1 \choose \nu_2}p_1^{\nu_1}p_2^{\nu_2}(1-p_1-p_2)^{n-\nu_1-\nu_2}\left.E\left(\left.\frac{-\partial^2\log
L(\theta)}{\partial\theta_r\partial\theta_s}\right|\nu_1,\nu_2\right)\right|_{\hat{\theta}}.
\end{eqnarray*}
One may write
\begin{eqnarray*}
\E\left(\frac{-\partial^2\log
L(\theta)}{\partial\theta_r\partial\theta_s}|\nu_1,\nu_2\right)&=&(\nu_1G_1(\tau)+\nu_2(G_2(c)-G_2(\tau)))\alpha_{r,s}+{\nu_1}\zeta_{1}(r,s)
+{\nu_2}\zeta_{2}(r,s)\\
&+&({n-\nu_1-\nu_2})\varphi(r,s),
\end{eqnarray*}
where
\[\alpha_{r,s}=\frac{-\partial^2}{\partial\theta_r\partial\theta_s}\log[(\sigma_X\sigma_Y)^{-1}(1-\rho^2)^{-1/2}],\]
\begin{eqnarray*}
\zeta_{j}(r,s)&=&\int_{-\infty}^{\infty}\int_{0}^{\infty} t^{-1}\frac{\partial^2Q_{j}(y,t)}{\partial\theta_r\partial\theta_s}
P_{f_j}(y,t;\theta)\;dt\;dy\\
&=&{\rm E}_{P_{f_j}}(T^{-1}\frac{\partial^2Q_{j}(Y,T)}{\partial\theta_r\partial\theta_s}),\quad j=1,2,\quad\mbox{say},
\end{eqnarray*}
\[\varphi(r,s)=\int_{-\infty}^{\infty}h(y;r,s)\;dy+\int_{-\infty}^{\infty}[g(y;r)g(y;s)]/[P_{C2}(y;\theta)]\;dy\]
and
\[h(y;r,s)=\frac{-\partial^2P_{C2}(y;\theta)}{\partial\theta_r\partial\theta_s},\quad g(y;r)=
\frac{\partial P_{C2}(y;\theta)}{\partial\theta_r}.\]
The functions $\alpha_{r,s}$ and $\zeta_{j}(r,s)$ for $j=1,2$ are simplified and given in the Appendix. It is straightforward that if
$T^{-1}\frac{\partial^2Q_{j}(Y,T)}{\partial\theta_r\partial\theta_s}$ is a function of $T$ only, the expectation can be taken on $f_T$
instead of $P_{f_j}$. The functions $h(y;r,s)$ and $g(y;r)$ are simplified as
\[g(y;r)=\sum_{k=1}^{2}\sum_{j_1=0}^{1}\sum_{j_2=1}^{2-j_1}(-1)^{k-1}e^{(k-1)\beta_2}
\lambda(y;r,k,j_1,j_2)\Phi^{(j_1)}(c_{2}(y;k,1))\Phi^{(j_2)}(c_{2}(y;k,2)),\]
and
\[h(y;r,s)=\sum_{k=1}^{2}\sum_{j_1=0}^{2}\sum_{j_2=1}^{3-j_1}(-1)^{k}e^{(k-1)\beta_2}
\gamma(y;r,s,k,j_1,j_2)\Phi^{(j_1)}(c_{2}(y;k,1))\Phi^{(j_2)}(c_{2}(y;k,2)),\]
where, $\Phi^{(j)}$ is the $j^{\mbox{th}}$ derivative of $\Phi$
and the coefficients $\lambda(y;r,k,j_1,j_2)$ and
$\gamma(y;r,s,k,j_1,j_2)$ are given in the Appendix.

\section{Illustrative example}

In order to illustrate the results of previous sections, let us study a numerical example.
Whitmore et al., 1998 presents a real data set on failure age and three potential markers for aluminum
reduction cells in a Canadian aluminum smelter. The production process of Aluminum consists of
electrolysis of molten alumina and cryolite in reduction cells.
Cryolite lowers the melting point of alumina to
$S_0=$950$^{\circ}$C. The cell's cathode is a carbon-lined steel box
which is subject to severe thermal, chemical and mechanical
stresses. The degradation of these cells can be marked by physical
distortion of the steel box. Suppose that $n=29$ reduction cells
are subjected to a step stress accelerated life test with $m=2$
stress levels $S_1=$1200$^{\circ}$C and $S_2=$1400$^{\circ}$C. Table
\ref{Alu2} provides an example of marker and failure data for 29 cells of a
particular design that were operated to failure under uniform
conditions in the Aluminum smelter. The censoring time is set to
$C=700$ days. The table shows the failure age (in days of service)
and the values at failure age of a marker for each cell, namely, the
horizontal distortion of the steel box (in inches). For these data
the threshold is taken to be $D=1$, the stress changing time is
$\tau=400$ days and no item is censored. We use these data to illustrate
the theoretical results of the optimization procedure.

\begin{table}
\caption{Data on failure age and horizontal distortion of the box as
a marker for aluminum reduction cells.}\label{Alu2}
\begin{tabular}{c | c c c c c c c c c c}
\hline \hline
Cell   &         1     & 2   &  3   &  4   &  5   &  6   &  7   &  8   &  9  & 10 \\
Stress level   & $S_2$     & $S_2$   & $S_1$   &  $S_2$   &  $S_2$   &  $S_1$   & $S_2$   &  $S_1$   & $S_1$  & $S_2$ \\
\hline
Failure Age     &573    &  447   &   365    &  412    &  508    &  385   &   611   &   235   &   395   &   471  \\
(in days) & & & & & & & & \\
Horiz. Distort. &  4.16    &   2.71   &    2.17   &    3.89   &    4.22   &    4.14   &    4.66   &    2.53   &    2.73   &    1.91 \\
(in inches) & & & & & & & & \\
\hline \hline
Cell   &        11   &   12   &  13  & 14   &  15  & 16   &  17  &  18 & 19 & 20  \\
Stress level   & $S_2$     & $S_2$   & $S_2$   &  $S_1$   &  $S_2$   &  $S_1$   & $S_2$   &  $S_2$   & $S_1$  & $S_2$ \\
\hline
Failure Age   & 604   &   509  &    653   &   341   &   441   &   392   &   447   &   486   &   341   &   666 \\
(in days) & & & & & & & & \\
Horiz. Distort. & 4.40  &    4.61   &   2.57   &   3.65   &   2.82   &   3.00   &   3.05   &   3.33   &   1.82   &   4.02\\
(in inches) & & & & & & & &\\
\hline \hline
Cell   &        21   & 22   &  23  & 24   &  25  & 26   &  27  &  28 & 29 &   \\
Stress level   & $S_2$     & $S_1$   & $S_2$   &  $S_2$   &  $S_2$   &  $S_2$   & $S_2$   &  $S_2$   & $S_2$  &  \\
\hline
Failure Age    & 589   &    347    &   588    &   577   &    567   &    468    &   564   &    435   &    504   &    \\
(in days) & & & & & & & & \\
Horiz. Distort. & 4.11    &   2.41   &    3.27   &    4.36    &   2.95    &   2.90   &    3.58    &   1.75   &    3.95     &   \\
(in inches) & & & & & & & &\\
\end{tabular}
\end{table}

\begin{figure}[!hbtp]
\centerline{\psfig{figure=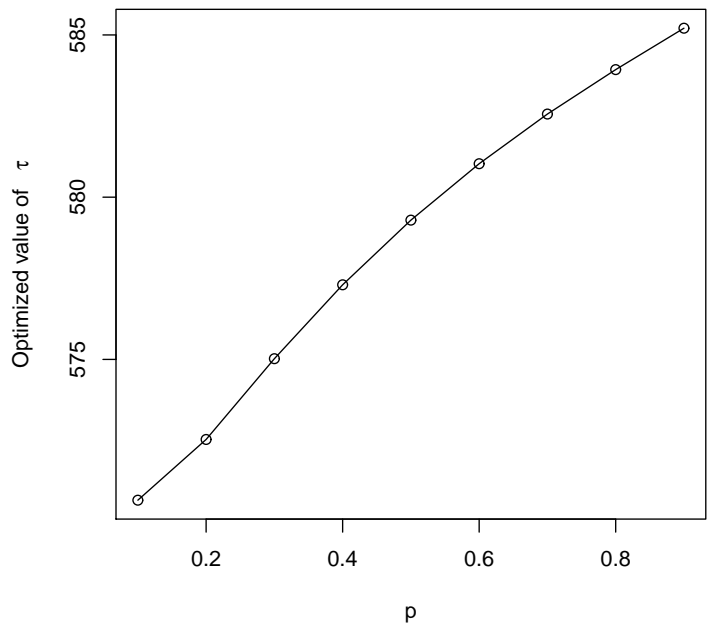,width=10cm}}
 \caption{\scriptsize Values of $\tau^*$ for different values of $p$.}\label{plot1}
\end{figure}

\begin{table}
\centering \caption{Optimal SSALT plan for minimizing
$\mbox{Avar}(\hat{\xi}_{p})$ for different values of $p$.\label{opt}}
\begin{tabular}{c | c c c c c }
\hline\hline
$p$ & $\hat{\xi}_{p}$ & minimum C.V. & $\tau^{*}$ & $G_1(\tau^{*})$ & $G_2(C)-G_2(\tau^{*})$\\
\hline
0.1 & 286.0 & 1.102 & 570.66 & 0.3197 & 0.0775\\
0.2 & 442.1 & 1.556 & 572.53 & 0.3208 & 0.0763\\
0.3 & 630.4 & 2.050 & 575.02 & 0.3222 & 0.0747\\
0.4 & 878.6 & 2.620 & 577.30 & 0.3235 & 0.0732\\
0.5 & 1227.6&3.304  & 579.29 & 0.3247 & 0.0719\\
0.6 & 1753.1& 4.155 & 581.03 & 0.3257 & 0.0708\\
0.7 & 2618.5& 5.257 & 582.56 & 0.3266 & 0.0698\\
0.8 & 4256.4& 6.775 & 583.93 & 0.3274 & 0.0689\\
0.9 & 8350.1&9.110  & 585.21 & 0.3281 & 0.0681\\
\hline
\end{tabular}
\end{table}

\begin{figure}[!hbtp]
\centerline{\psfig{figure=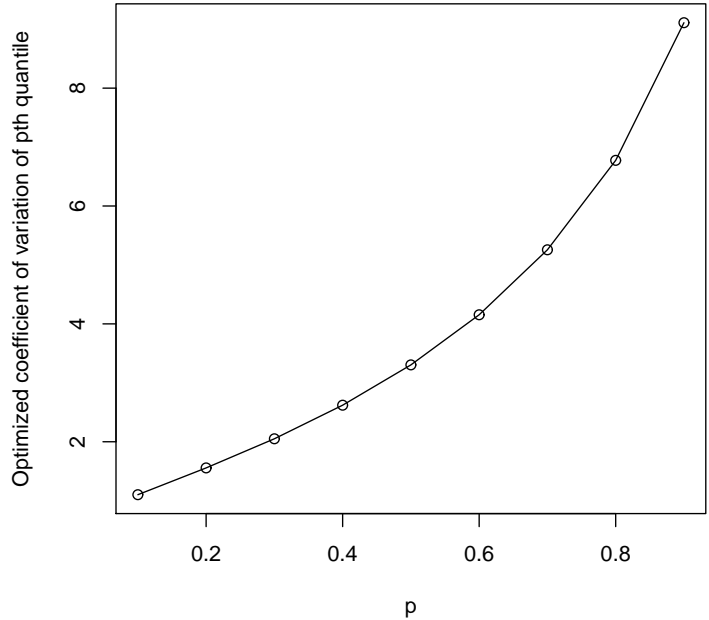,width=10cm}}
 \caption{\scriptsize The optimized approximated coefficient of variation of $\hat{\xi}_p$ for different values of $p$.}\label{plot2}
\end{figure}

Using these data, one can obtain the ML estimates of the parameters
using the likelihood in \eqref{mainloglike} as $\hat{\mu}_{X1}=0.0005$,
$\hat{\mu}_{X2}=0.0007$, $\hat{\mu}_{Y1}=0.0005$, $\hat{\mu}_{Y2}=0.0006$,
$\hat{\sigma}_X=0.0011$, $\hat{\sigma}_Y=0.0188$ and $\hat{\rho}=0.9422$.

The optimization process for minimizing $\mbox{Avar}(\hat{\xi}_{p})$
is performed using the optimization procedures of software R.2.14.1.
The results including $\hat{\xi}_{p}$, the optimized approximated
coefficient of variation of $\hat{\xi}_{p}$ (minimum C.V.), the
optimized time $\tau^{*}$, the probability of failure under the
stress level $S_1$ that is $G_1(\tau^{*})$ and the probability of
failure under the stress level $S_2$, i.e. $G_2(C)-G_2(\tau^{*})$,
are obtained for $p=0.1(0.1)0.9$ and tabulated in Table \ref{opt}.

Figure \ref{plot1} shows the plot of $\tau^*$ as a function of $p$. The values of
the optimized approximated coefficient of variation of $\hat{\xi}_{p}$ are also plotted for different values of
$p$ in Figure \ref{plot2}. As it can be seen from Figures \ref{plot1} and \ref{plot2}, the optimal time $\tau^*$ is an increasing function
of $p$. It is legal to have such a result, since under a higher stress level the items fail more rapidly and such failures contain
more information about lower quantiles of the lifetime distribution of the products. As one can observe from Figure \ref{plot2},
the precision of the optimal estimate of $\hat{\xi}_{p}$ decreases for the upper percentiles of the products' lifetime distribution.

\section*{Acknowledgements}

The authors would like to thank two anonymous referees and also
Professor M. N. Tata from Shahid Bahonar University and
Professor Massoud Amini from Tarbiat Modares University for their
valuable comments and suggestion which improved the contents of this
paper.

\section*{Appendix A. Fisher information matrix}

Denoting $\frac{\partial f}{\partial\theta_r}$ and $\frac{\partial^2
f}{\partial\theta_r\partial\theta_s}$ by $f^{[r]}$ and $f^{[r,s]}$,
respectively, we have
\[\lambda(y;r,1,0,1)=c_y^{[r]},\quad \lambda(y;r,1,1,1)=c_y c_{2}^{[r]}(y;1,1),
\quad \lambda(y;r,1,0,2)=c_y c_{2}^{[r]}(y;1,2),\]\vspace{-4mm}
\[\lambda(y;r,2,0,1)=c_y^{[r]}+\beta_2^{[r]}c_y,\quad \lambda(y;r,2,1,1)=c_y c_{2}^{[r]}(y;2,1),
\quad \lambda(y;r,2,0,2)=c_y c_{2}^{[r]}(y;2,2),\]\vspace{-4mm}
\[\gamma(y;r,s,1,0,1)=c_y^{[r,s]},\quad \gamma(y;r,s,1,1,1)=c_y c_{2}^{[r,s]}(1,1)+
c_y^{[r]} c_{2}^{[s]}(y;1,1)+c_y^{[s]}
c_{2}^{[r]}(y;1,1),\]\vspace{-4mm}
\[\gamma(y;r,s,1,0,2)=c_y c_{2}^{[r,s]}(y;1,2)+
c_y^{[r]} c_{2}^{[s]}(y;1,2)+c_y^{[s]}
c_{2}^{[r]}(y;1,2),\]\vspace{-4mm}
\[\gamma(y;r,s,k,1,2)=c_y[ c_{2}^{[r]}(y;k,1)c_{2}^{[s]}(y;k,2)+c_{2}^{[s]}(y;k,1) c_{2}^{[r]}(y;k,2)],\quad k=1,2,\]\vspace{-4mm}
\[\gamma(y;r,s,k,2,1)=c_y c_{2}^{[r]}(y;k,1)c_{2}^{[s]}(y;k,1),\quad
\gamma(y;r,s,k,0,3)=c_y c_{2}^{[r]}(y;k,2)c_{2}^{[s]}(y;k,2),,\quad
k=1,2,\]\vspace{-4mm}
\[\gamma(y;r,s,2,0,1)=c_y^{[r,s]}+\beta_2^{[s]}c_y^{[r]}+\beta_2^{[r]}c_y^{[s]}+c_y\beta_2^{[r]}\beta_2^{[s]}
+c_y\beta_2^{[r,s]},\]\vspace{-4mm}
\[\gamma(y;r,s,2,1,1)=\beta_2^{[s]}c_{2}^{[r]}(y;2,1)c_y+c_{2}^{[r,s]}(y;2,1)c_y+c_{2}^{[r]}(y;2,1)c_y^{[s]}+
c_{2}^{[s]}(y;2,1)c_y^{[r]}+\beta_2^{[r]}c_{2}^{[s]}(y;2,1)c_y,\]\vspace{-4mm}
\[\gamma(y;r,s,2,0,2)=\beta_2^{[s]}c_{2}^{[r]}(y;2,2)c_y+c_{2}^{[r,s]}(y;2,2)c_y+c_{2}^{[r]}(y;2,2)c_y^{[s]}+
c_{2}^{[s]}(y;2,2)c_y^{[r]}+\beta_2^{[r]}c_{2}^{[s]}(y;2,2)c_y,\]
Letting $\eta_4=\rho(1-\rho^2)^{-1},\;\eta_5=\rho^{-1}+\eta_4$ and
$\eta_6=\rho^{-1}+2\eta_4,$ we have
\[c_y^{[6]}=-\sigma_Y^{-2}c_y/2,\;c_y^{[6,6]}=3\sigma_Y^{-4}c_y/4,\;\;c_y^{[r,s]}=c_y^{[r]}=0,\quad\mbox{for }r\neq6,\;s\ne6,\]
\[c_{2}^{[1]}(y;1,1)=c_{2}^{[1]}(y;2,1)=-\tau\eta_3,\;c_{2}^{[2]}(y;1,1)=c_{2}^{[2]}(y;2,1)=-(C-\tau)\eta_3,\]\vspace{-4mm}
\[c_{2}^{[3]}(y;1,1)=c_{2}^{[3]}(y;2,1)=\tau\eta_2\eta_3,\;c_{2}^{[4]}(y;1,1)=c_{2}^{[4]}(y;2,1)=(C-\tau)\eta_2\eta_3,\]\vspace{-4mm}
\[c_{2}^{[5]}(y;1,1)=\eta_3\sigma_X^{-2}[\eta_2q_{2}(C,y)-P_{2}(C)/2],\;
c_{2}^{[5]}(y;2,1)=\eta_3\sigma_X^{-2}[\eta_2q_{2}(C,y)-(P_{2}(C)-2D(1-\rho^2))/2],\]\vspace{-4mm}
\[c_{2}^{[6]}(y;1,1)=c_{2}^{[6]}(y;2,1)=-\sigma_Y^{-2}\eta_2\eta_3q_{2}(C,y)/2,\]\vspace{-4mm}
\[c_{2}^{[7]}(y;1,1)=\eta_3[\eta_4P_{2}(C)-\eta_2\eta_5q_{2}(C,y)],\;
c_{2}^{[7]}(y;2,1)=\eta_3[\eta_4P_{2}(C)-\eta_2\eta_5q_{2}(C,y)+2D\rho],\]\vspace{-4mm}
\[c_{2}^{[1,5]}(y;1,1)=c_{2}^{[1,5]}(y;2,1)=\frac{\tau}{2}\sigma^{-2}_X\eta_3,\;
c_{2}^{[1,7]}(y;1,1)=c_{2}^{[1,7]}(y;2,1)=-\tau\eta_3\eta_4,\]\vspace{-4mm}
\[c_{2}^{[2,5]}(y;1,1)=c_{2}^{[2,5]}(y;2,1)=\frac{(C-\tau)}{2}\sigma^{-2}_X\eta_3,\;
c_{2}^{[2,7]}(y;1,1)=c_{2}^{[2,7]}(y;2,1)=-(C-\tau)\eta_3\eta_4,\]\vspace{-4mm}
\[c_{2}^{[3,5]}(y;1,1)=c_{2}^{[3,5]}(y;2,1)=-\tau\sigma^{-2}_X\eta_3\eta_2,\;
c_{2}^{[3,6]}(y;1,1)=c_{2}^{[3,6]}(y;2,1)=\frac{\tau}{2}\sigma^{-2}_Y\eta_3\eta_2,\]\vspace{-4mm}
\[c_{2}^{[3,7]}(y;1,1)=c_{2}^{[3,7]}(y;2,1)=\tau\eta_3\eta_2\eta_5,\;
c_{2}^{[4,5]}(y;1,1)=c_{2}^{[4,5]}(y;2,1)=-(C-\tau)\sigma^{-2}_X\eta_3\eta_2,\]\vspace{-4mm}
\[c_{2}^{[4,6]}(y;1,1)=c_{2}^{[4,6]}(y;2,1)=\frac{(C-\tau)}{2}\sigma^{-2}_Y\eta_3\eta_2,\;
c_{2}^{[4,7]}(y;1,1)=c_{2}^{[4,7]}(y;2,1)=(C-\tau)\eta_3\eta_2\eta_5,\]\vspace{-4mm}
\[c_{2}^{[5,5]}(y;1,1)=\eta_3\sigma^{-4}_X[\frac{3}{4}P_2(C)-2\eta_2q_{2}(C,y)],\;
c_{2}^{[5,5]}(y;2,1)=\eta_3\sigma^{-4}_X[\frac{3}{4}(P_2(C)-2D(1-\rho^2))-2\eta_2q_{2}(C,y)],\]\vspace{-4mm}
\[c_{2}^{[5,6]}(y;1,1)=c_{2}^{[5,6]}(y;2,1)=\frac{1}{2}\sigma^{-2}_X\sigma_Y^{-2}\eta_3\eta_2q_{2}(C,y),\;
c_{2}^{[5,7]}(y;1,1)=\sigma_X^{-2}\eta_3[\eta_2\eta_5q_{2}(C,y)-\frac{1}{2}\eta_4P_2(C)],\]\vspace{-4mm}
\[c_{2}^{[5,7]}(y;2,1)=\sigma_X^{-2}\eta_3[\eta_2\eta_5q_{2}(C,y)-\frac{1}{2}\eta_4P_2(C)-\rho D],\;
c_{2}^{[6,6]}(y;1,1)=c_{2}^{[6,6]}(y;2,1)=\frac{1}{4}\sigma_Y^{-4}\eta_3\eta_2q_{2}(C,y),\]\vspace{-4mm}
\[c_{2}^{[6,7]}(y;1,1)=c_{2}^{[6,7]}(y;2,1)=-\frac{1}{2}\sigma_Y^{-2}\eta_3\eta_2\eta_5q_{2}(C,y),\]\vspace{-4mm}
\[c_{2}^{[7,7]}(y;1,1)=\eta_3[(1+2\rho^2)(1-\rho^2)^{-2}P_2(C)-3\eta_2(1-\rho^2)^{-2}q_{2}(C,y)],\]\vspace{-4mm}
\[c_{2}^{[7,7]}(y;2,1)=\eta_3[(1+2\rho^2)(1-\rho^2)^{-2}P_2(C)-3\eta_2(1-\rho^2)^{-2}q_{2}(C,y)+2D(1-\rho^2)^{-1}],\]
and $c_{2}^{[r,s]}(y;1,1)=c_{2}^{[r,s]}(y;2,1)=0$ otherwise $r,s$.
Also,
\[c_{2}^{[3]}(y;1,2)=c_{2}^{[3]}(y;2,2)=-\tau c_y,\;
c_{2}^{[4]}(y;1,2)=c_{2}^{[4]}(y;2,2)=-(C-\tau)c_y,\]
\[c_{2}^{[5]}(y;1,2)=0,\; c_{2}^{[5]}(y;2,2)=D\rho
C^{-1/2}\sigma_X^{-3},\]\vspace{-4mm}
\[c_{2}^{[6]}(y;1,2)=c_{2}^{[6]}(y;2,2)=-\frac{1}{2}\sigma_Y^{-2}q_{2}(C,y)c_y,\;
c_{2}^{[7]}(y;1,2)=0,\;
c_{2}^{[7]}(y;2,2)=-2DC^{-1/2}\sigma_X^{-1},\]\vspace{-4mm}
\[c_{2}^{[r]}(y;1,2)=c_{2}^{[r]}(y;2,2)=0,\; r=1,2,\]\vspace{-4mm}
\[c_{2}^{[3,6]}(y;1,2)=c_{2}^{[3,6]}(y;2,2)=\frac{\tau}{2}\sigma^{-2}_Yc_y,\;
c_{2}^{[4,6]}(y;1,2)=c_{2}^{[4,6]}(y;2,2)=\frac{(C-\tau)}{2}\sigma^{-2}_Yc_y,\]\vspace{-4mm}
\[c_{2}^{[5,5]}(y;1,2)=0,\;
c_{2}^{[5,5]}(y;2,2)=-\frac{3}{2}D\sigma_X^{-5}C^{-1/2}\rho,\;
c_{2}^{[5,7]}(y;1,2)=0,\;
c_{2}^{[5,7]}(y;2,2)=D\sigma_X^{-3}C^{-1/2},\]\vspace{-4mm}
\[c_{2}^{[6,6]}(y;1,2)=c_{2}^{[6,6]}(y;2,2)=\frac{3}{4}\sigma_Y^{-4}q_{2}(C,y)c_y,\]
and $c_{2}^{[r,s]}(y;1,2)=c_{2}^{[r,s]}(y;2,2)=0$ otherwise $r,s$.
Furthermore
\[\beta_2^{[1]}=2D\tau\sigma_X^{-2}C^{-1},\;\beta_2^{[2]}=2D(C-\tau)\sigma_X^{-2}C^{-1},\;
\beta_2^{[5]}=-2D(D-P_2(C))\sigma_X^{-4}C^{-1},\]\vspace{-4mm}
\[\beta_2^{[r]}=0,\;r=3,4,6,7,\]\vspace{-4mm}
\[\beta_2^{[1,5]}=-2D\tau\sigma_X^{-4}C^{-1},\;\beta_2^{[2,5]}=-2D(C-\tau)\sigma_X^{-4}C^{-1},\;
\beta_2^{[5,5]}=4D(D-P_2(C))\sigma_X^{-6}C^{-1},\] and
$\beta_2^{[r,s]}=0,$ otherwise $r,s$.
\[\alpha_{5,5}=-\frac{1}{2\sigma_X^4},\;
\alpha_{6,6}=-\frac{1}{2\sigma_Y^4},\;
\alpha_{7,7}=-\frac{(1+\rho^2)}{(1-\rho^2)^2},\; \alpha_{r,s}=0\]
otherwise $r,s$.
 Also, for $j=1,2,$
\[\zeta_{j}(1,1)=(1-\rho^2)^{-1}\sigma_X^{-2}{\rm E}(T^{3-2j})\tau^{2j-2},\;
\zeta_{1}(1,2)=0,\;
\zeta_{2}(1,2)=(1-\rho^2)^{-1}\sigma_X^{-2}\tau(1-\tau {\rm E}(T^{-1})),\]\vspace{-4mm}
\[\zeta_{j}(1,3)=-2\eta_1\eta_2{\rm E}(T^{3-2j})\tau^{2j-2},\;
,\zeta_{1}(1,4)=0,\;
\zeta_{2}(1,4)=-2\eta_1\eta_2\tau(1-\tau {\rm E}(T^{-1})),\]\vspace{-4mm}
\[\zeta_{j}(1,5)=\frac{1}{2}\sigma_X^{-4}(1-\rho^2)^{-1}[2{\rm E}(T^{1-j}P_j(T))-\rho\sigma_X\sigma_Y^{-1}{\rm E}_{P_{f_j}}(T^{1-j}q_{j}(T,Y))]\tau^{j-1},\]
\[\zeta_{j}(1,6)=-\sigma_Y^{-2}\eta_1\eta_2{\rm E}_{P_{f_j}}(T^{1-j}q_{j}(T,Y))\tau^{j-1},\]\vspace{-4mm}
\[\zeta_{j}(1,7)=2\eta_2\eta_1[\eta_6{\rm E}_{P_{f_j}}(T^{1-j}q_{j}(T,Y))-2\eta_2\eta_4{\rm E}(T^{1-j}P_j(T))]\tau^{j-1},\]
\[\zeta_{1}(2,2)=\zeta_{1}(2,3)=\cdots=\zeta_{1}(2,7)=0,\]\vspace{-4mm}
\[\zeta_{2}(2,2)=(1-\rho^2)^{-1}\sigma_X^{-2}{\rm E}(T^{-1}(T-\tau)^2),\;\zeta_{2}(2,3)=-2\eta_1\eta_2\tau(1-\tau {\rm E}(T^{-1})),\]\vspace{-4mm}
\[\zeta_{2}(2,4)=-2\eta_1\eta_2{\rm E}(T^{-1}(T-\tau)^2),\]
\begin{eqnarray*}
\zeta_{2}(2,5)&=&\frac{1}{2}\sigma_X^{-4}(1-\rho^2)^{-1}[2({\rm E}(P_2(T))-\tau {\rm E}(T^{-1}P_2(T)))\\
&&-\rho\sigma_X\sigma_Y^{-1}({\rm E}_{P_{f_2}}(q_2(T,Y))-\tau {\rm E}_{P_{f_2}}(T^{-1}q_2(T,Y)))],
\end{eqnarray*}\vspace{-4mm}
\[\zeta_{2}(2,6)=-\sigma_Y^{-2}\eta_1\eta_2({\rm E}_{P_{f_2}}(q_2(T,Y))-\tau {\rm E}_{P_{f_2}}(T^{-1}q_2(T,Y))),\]\vspace{-4mm}
\begin{eqnarray*}
\zeta_{2}(2,7)&=&2\eta_2\eta_1[\eta_6({\rm E}_{P_{f_j}}(q_j(T,Y))-\tau {\rm E}_{P_{f_j}}(T^{-1}q_j(T,Y)))\\
&&-2\eta_2\eta_4({\rm E}(P_j(T))-\tau {\rm E}(T^{-1}P_j(T)))],
\end{eqnarray*}\vspace{-4mm}
\[\zeta_{j}(3,3)=2\eta_1{\rm E}(T^{3-2j})\tau^{2j-2},\]\vspace{-4mm}
\[\zeta_{1}(3,4)=0,\;\zeta_{2}(3,4)=2\eta_1\tau(1-\tau {\rm E}(T^{-1})),\;
\zeta_{j}(3,5)=-\sigma_X^{-2}\eta_2\eta_1{\rm E}(T^{1-j}P_j(T))\tau^{j-1},\]\vspace{-4mm}
\[\zeta_{j}(3,6)=\sigma_Y^{-2}\eta_1[2{\rm E}_{P_{f_j}}(T^{1-j}q_{j}(T,Y))-\eta_2{\rm E}(T^{1-j}P_j(T))]\tau^{j-1},\]
\[\zeta_{j}(3,7)=2\eta_1[\eta_2\eta_6{\rm E}(T^{1-j}P_j(T))-2\eta_4{\rm E}_{P_{f_j}}(T^{1-j}q_{j}(T,Y))]\tau^{j-1},\]\vspace{-4mm}
\[\zeta_{1}(4,4)=\zeta_{1}(4,5)=\zeta_{1}(4,6)=\zeta_{1}(4,7)=0,\;\zeta_{2}(4,4)=2\eta_1{\rm E}(T^{-1}(T-\tau)^2),\]\vspace{-4mm}
\[\zeta_{2}(4,5)=-\sigma_X^{-2}\eta_2\eta_1({\rm E}(P_2(T))-\tau {\rm E}(T^{-1}P_2(T))),\]
\[\zeta_{2}(4,6)=\sigma_Y^{-2}\eta_1[2({\rm E}_{P_{f_2}}(q_2(T,Y))-\tau {\rm E}_{P_{f_2}}(T^{-1}q_2(T,Y)))-\eta_2({\rm E}(P_2(T))-\tau {\rm E}(T^{-1}P_2(T)))],\]\vspace{-4mm}
\[\zeta_{2}(4,7)=2\eta_1[\eta_2\eta_6({\rm E}(P_2(T))-\tau {\rm E}(T^{-1}P_2(T)))-2\eta_4({\rm E}_{P_{f_2}}(q_2(T,Y))-\tau {\rm E}_{P_{f_2}}(T^{-1}q_2(T,Y)))],\]\vspace{-4mm}
\[\zeta_{j}(5,5)=\sigma_X^{-4}\eta_1\eta_2[2\eta_2{\rm E}(T^{-1}P_j^2(T))-3{\rm E}_{P_{f_j}}(T^{-1}P_j(T)q_{j}(T,Y)))/2]+\sigma_X^{-6}{\rm E}(T^{-1}P_j^2(T)),\]\vspace{-4mm}
\[\zeta_{j}(5,6)=\frac{-1}{2}\sigma_X^{-2}\sigma_Y^{-2}\eta_1\eta_2{\rm E}_{P_{f_j}}(T^{-1}P_j(T)q_{j}(T,Y))),\]
\[\zeta_{j}(5,7)={\sigma_X^{-2}}\eta_1\eta_2[\eta_6{\rm E}_{P_{f_j}}(T^{-1}P_j(T)q_{j}(T,Y)))-2\eta_5\eta_2{\rm E}(T^{-1}P_j^2(T))],\]\vspace{-4mm}
\[\zeta_{j}(6,6)=\eta_1\sigma_Y^{-4}[2{\rm E}_{P_{f_j}}(T^{-1}q_{j}^2(T,Y)))-\frac{3}{2}\eta_2{\rm E}_{P_{f_j}}(T^{-1}P_j(T)q_{j}(T,Y)))],\]
\[\zeta_{j}(6,7)=-\eta_1\sigma_Y^{-2}[2\eta_4{\rm E}_{P_{f_j}}(T^{-1}q_{j}^2(T,Y)))-\eta_2\eta_6{\rm E}_{P_{f_j}}(T^{-1}P_j(T)q_{j}(T,Y)))],\]\vspace{-4mm}
\begin{eqnarray*}
\zeta_{j}(7,7)&=&2\eta_1\eta_4^2((3+\rho^{-2}){\rm E}_{P_{f_j}}(T^{-1}q_{j}^2(T,Y)))+(3\rho^{-2}+\rho^{-4})(\eta_2{\rm E}(T^{-1}P_j^2(T)))\\
&&-2(1+3\rho^{-2})\eta_2{\rm E}_{P_{f_j}}(T^{-1}P_j(T)q_{j}(T,Y)))).
\end{eqnarray*}\vspace{-4mm}

\begin{thebibliography}{99}
\bibitem{bn}{Bagdonavicius, V. and Nikulin, M. (2010).} {\itshape Accelerated Life Models: Modeling and Statistical Analysis}. Taylor \& Francis.

\bibitem{chf}{Chhikara R.S. and Folks, J.L. (1989).} {\itshape The Inverse Gaussian Distribution: Theory, Methodology and Applocations}.
Marcel Dekker: New York.


\bibitem{cm}{Cox D.R. and Miller, H.D. (1965).} {\itshape The Theory of Stochastic Processes}. John Wiley \& Sons: New York.

\bibitem{jk}{Jewell, N.P. and Kalbfleisch, J.D. (1996).} {Marker Processes in Survival Analysis}. {\em Lifetime Data Analysis} { 2}:15--29.

\bibitem{jin} Jin, M. (2011). Analysis of Failure Time Data with Mixed-Effects Accelerated Failure Time Model.
{\em Communications in Statistics -- Simulation and Computation} {40}:614--619.

\bibitem{lt}{Liao, C.M. and Tseng, C.T. (2006).} {Optimal Design for Step-Stress Accelerated Degradation Tests}.
{\em IEEE Transactions on Reliability} {55}:59--66.

\bibitem{lu}{Lu, J. (1995).} {\itshape A Reliability Model Based on Degradation and Lifetime Data}. Ph.D. Thesis, McGill University, Montreal, Canada.

\bibitem{nel}{Nelson, W. (1990).} {\itshape Accelerated Testing: Statistical Models, Test Plans, and Data Analysis}. Wiley: New York.

\bibitem{pb}{Pan, Z. and Balakrishnan, N. (2010).} {Multiple Steps Step-Stress Accelerated Degradation Modeling Baesd on Wiener and Gamma Process}.
{\em Communication in Statistics-Simulation and Computation} {39}:1384--1402.

\bibitem{pbs} Pan, Z., Balakrishnan, N. and Sun, Q. (2011). Bivariate Constant-Stress Accelerated Degradation Model and Inference. {\em Communications in Statistics - Simulation and Computation} {40}:247--257.

\bibitem{shm} Simino, J., Hollander, M. and McGee, D. (2012). Calibration of Proportional Hazards and Accelerated Failure Time Models.
{\em Communications in Statistics -- Simulation and Computation} {41}:922--941.

\bibitem{tyx}{Tang, L.C., Yang, G.Y. and Xie, M. (2004).} {Planning of step-stress accelerated degradation test}. {\em 2004 Annual Symposium - RAMS}. 287–-292.

\bibitem{tbt}{Tseng, S.T., Balakrishnan, N. and Tsai, C.C. (2009).} {Optimal Step-Stress Accelerated Degradation Test Plan for
Gamma Degradation Process}. {\em IEEE Transactions on Reliability } {58}:611--618.

\bibitem{wcl} Wang, F. K., Cheng, Y. F. and Lu, W. L. (2012). Partially Accelerated Life Tests for the Weibull Distribution Under Multiply Censored Data. {\em    Communications in Statistics - Simulation and Computation} {41}:1667--1678.

\bibitem{wcl} {Whitmore, G.A. Crowder, M.J. and Lawless, J.F. (1998).} {Failure Inference from a Marker Process based on Bivariate Wiener Model},
{\em Lifetime Data Analysis} {4}:229--251.

\bibitem{my}{Yashin, A.I. and Manton, K.G. (1997).} {Effects of unobserved and partially observed covariate processes on system failure:
 a review of models and estimation strategies}. {\em Statistical Science} {12}:20--34.
\end{thebibliography}
\end{document}